\documentclass[leqno,11pt,a4paper,twoside]{article}%
\usepackage{amssymb}
\usepackage{amsfonts}
\usepackage{amsmath}
\usepackage{authblk}
\usepackage{graphicx}
\usepackage{ifthen}
\usepackage{comment}
\usepackage{xcolor}
\usepackage{enumitem}
\usepackage{lipsum,multicol}
\usepackage[T1]{fontenc}
\usepackage[english]{babel}
\usepackage[english]{babel}
\setlist[itemize]{noitemsep, topsep=0pt}
\setlist[enumerate]{noitemsep, topsep=0pt}
\setlist[itemize]{leftmargin=*}
\setlist[enumerate]{leftmargin=*}
\providecommand{\U}[1]{\protect\rule{.1in}{.1in}}

\providecommand{\norm}[1]{\left\lVert#1\right\rVert}
\providecommand{\abs}[1]{\left\lvert#1\right\rvert}
\providecommand{\pr}[1]{\left(#1\right)} 
\providecommand{\pp}[1]{\left[#1\right]} 
\providecommand{\set}[1]{\left\lbrace#1\right\rbrace} 
\providecommand{\scal}[1]{\left\langle#1\right\rangle}

\newcommand{\Lo}[2]{\mathbb{L}^{#1}\pr{{\mathcal{O}_{{#2}}}}}
\newcommand{\Ho}[1]{H_0^1\pr{\mathcal{O}_{#1}}}
\newcommand{\Hmo}[1]{H^{-1}\pr{\mathcal{O}_{#1}}}
\newcommand{\normi}[3]{\norm{#1}_{
    \ifthenelse{\equal{#2}{1}}{H_0^1\pr{\mathcal{O}_{#3}}}{%
    \ifthenelse{\equal{#2}{-1}}{H^{-1}\pr{\mathcal{O}_{#3}}}{}}}}
    \newcommand{\X}{X^{t,\xi,\eta,u}}
    \newcommand{\Y}{Y^{t,\xi,\eta,u}}
    \newcommand{\Xp}{X^{t,\xi',\eta',u}}
    \newcommand{\Yp}{Y^{t,\xi',\eta',u}}
    \newcommand{\XX}{\mathbb{X}^{t,\xi,\eta,u}}
     \newcommand{\YY}{\mathbb{Y}^{t,\xi,\eta,u}}
     \newcommand{\XXp}{\mathbb{X}^{t,\xi',\eta',u}}
     \newcommand{\YYp}{\mathbb{Y}^{t,\xi',\eta',u}}
      \newcommand{\n}[1]{\mid\norm{#1}\mid}
      \newcommand{\phiK}[2]{\underset{1\leq k\leq j}{\sum}\scal{#1,e_k^1}_{\Hmo{1}}^2-#2}

\makeatletter
\newcommand{\subjclass}[2][2020]{%
  \let\@oldtitle\@title%
  \gdef\@title{\@oldtitle\footnotetext{\textbf{#1 \emph{Mathematics subject classification.}} #2}}%
}
\newcommand{\keywords}[1]{%
  \let\@@oldtitle\@title%
  \gdef\@title{\@@oldtitle\footnotetext{\textbf{\emph{Key words.}} #1.}}%
}
\makeatother

\newcommand{\ef}[2]{\mathbb{E}\pp{#1\mid \mathcal{F}_{#2}}}

\newtheorem{theorem}{Theorem}

\newtheorem{corollary}[theorem]{Corollary}

\newtheorem{definition}[theorem]{Definition}

\newtheorem{proposition}[theorem]{Proposition}
\newtheorem{remark}[theorem]{Remark}

\newenvironment{proof}[1][Proof]{\noindent\textbf{#1.} }{\ \rule{0.5em}{0.5em}}

\begin{document}
\title{State-constrained porous media control systems with application to stabilization}
\author[1,2]{Ioana Ciotir}
\author[3,4]{Dan Goreac\footnote{Corresponding author, email: dan.goreac@univ-eiffel.fr}}
\author[5,6]{Ionu\c t Munteanu}
\affil[1]{Normandie University, INSA de Rouen Normandie,
LMI (EA 3226 - FR CNRS 3335), 76000 Rouen, France} 
\affil[2]{Research Center for Pure and Applied Mathematics,
Graduate School of Information Sciences,
Tohoku University, Japan}
\affil[3]{School of Mathematics and Statistics, Shandong University, Weihai, Weihai 264209, PR China} 
\affil[4]{LAMA, Univ Gustave Eiffel, UPEM, Univ Paris Est Creteil, CNRS, F-77447 Marne-la-Vallée, France}
\affil[5]{Faculty of Mathematics, Al. I. Cuza University, Bd. Carol I, 11, Iasi 700506, Romania}
\affil[6]{O. Mayer Institute of Mathematics, Romanian Academy, Bd. Carol I, 8, Iasi 700505, Romania}

\date{}


\maketitle

\begin{abstract}
{We aim at providing a characterization of the ability to maintain a stochastic coupled system with porous media components in a prescribed set of constraints by using internal controls.  This property is proven via a \textit{quasi-tangency} local-in-time condition in the spirit of Euler approximation schemes. In particular, by employing one of the components of the system as asymptotic supervisor, we give conditions guaranteeing the exponential asymptotic stabilizability of controlled porous media equations.}
\end{abstract}

\section{Introduction}

The present work focuses on the following porous media-type stochastic partial differential system with internal control
\begin{align}\label{Eq00}
\begin{cases}
dX^{t,\xi,\eta,u}(s)&=\pr{\Delta\beta_1\pr{X^{t,\xi,\eta,u}(s)}+f_1\pr{X^{t,\xi,\eta,u}(s),Y^{t,\xi,\eta,u}(s),u(s)}}ds\\[4pt]
&+\sigma_1\pr{X^{t,\xi,\eta,u}(s),Y^{t,\xi,\eta,u}(s),u(s)}dW(s),\\[4pt]
dY^{t,\xi,\eta,u}(s)&=\pr{\Delta\beta_2\pr{Y^{t,\xi,\eta,u}(s)}+f_2\pr{X^{t,\xi,\eta,u}(s),Y^{t,\xi,\eta,u}(s),u(s)}}ds\\[4pt] &+\sigma_2\pr{X^{t,\xi,\eta,u}(s),Y^{t,\xi,\eta,u}(s),u(s)}dW(s),\ s\geq t;\\[4pt]
X^{t,\xi,\eta,u}(t)&=\xi\in\mathbb{L}^2\pr{\Omega;\Hmo{1}},\ Y^{t,\xi,\eta,u}(t)=\eta\in\mathbb{L}^2\pr{\Omega;\Hmo{2}}.
\end{cases}
\end{align}
Here, for $i=1,2$, $\mathcal{O}_i$ are open bounded domains in euclidean spaces $\mathbb{R}^{d_i}$, $\beta_{i}$, for $i=1,2$ are Lipschitz monotone functions on $\mathbb{R}$,  and $W$ is a cylindrical Wiener process defined on some complete probability space $(\Omega, \mathcal{F}, (\mathcal{F}_{t})_{t\geq 0}, \mathbb{P})$, satisfying the usual conditions of completeness and right-continuity. The drift $f$ and the noise $\sigma$ coefficients satisfy suitable continuity properties to be made precise later on.

The stochastic porous media equation has been intensively studied recently in different frameworks. There is a vast literature on the subject and a complete overview would exceed the framework of this paper. We shall recall the pioneer works from \cite{DPRRWF}, \cite{Ren2006StochasticGP}, the case of an unbounded domain from \cite{BARBU20151024}, the critical cases from \cite{BabruLN}, \cite{CIOTIR2017595} and some properties of the solution from \cite{positivity}, \cite{BarbuExt2008}, \cite{GessExt2013}.
For further details, we invite our readers to refer to \cite{BDPR_2016} for a monograph on the subject. 

According to the order of growth of the operator $\beta$, the porous media equation describes different phenomena going from slow diffusion, for over-unit orders, to fast diffusion, for sub-unit orders, and even super-fast diffusion for negative orders. In the present work we are interested in the case when $\beta$ is monotonically increasing and Lipschitz continuous. The main interest of this case is that general maximal monotone graphs can be approximated by their Yosida approximation which are Lipschitz continuous and monotonically increasing. 
There are also significant physical problems with $\beta$ Lipschitz continuous as it is the case for the Stefan two phases problem.\\

The equations are controlled by a process $u$ and we give a \emph{state-constrained (or viability)} type result for the system. The interest of such questions for porous media type equations comes from geophysics and more precisely could be applied in the study and prevention of cliffs erosion. In other words, if we keep in mind that the solution of a porous media equation represents the moisture (volume of water over volume of voids), then the control can be seen as a drainage system which regulates this moisture. The system in our case could be seen as describing a cliff with two different types of structure corresponding to the two equations. We also consider that there exists a direct connection between the moisture and the solidity of the cliff. In this context, our type of result can be seen as a control methodology to reduce the cliff erosion.\\

{The literature on state-constrained systems originating from the seminal paper \cite{Nagumo_42} has known important advances in the past 30 years with the development of set-valued analysis (cf.  \cite{AubiCell84, aubin_frankowska_90}) and the extension of tangency concepts to stochastic finite-dimensional systems by \cite{AubinDaPrato_90}.  Two main streams exist relying either on tangency concepts in the Bouligand-Severi sense, e.g. \cite{aubin_frankowska_90, aubin_91, AubinDaPrato_90, Gautier_Thibault_93, Ciotir_Rascanu_2009} or on viscosity solutions, e.g. in \cite{BPQR1998, bardi_jensen_2002}.  Adapting the notion of tangency, the results in \cite{CarjaNeculaVrabie2007} offer an important semi-group-based method to deal with a wide class of deterministic PDEs or other infinite-dimensional components (e.g. mean-field control in \cite{BGL_2022}).  It is this later point of view we adopt here. 

The notion of quasi-tangency (see Definition \ref{DefLambdaQT}) provides a local quasi-constrained $\varepsilon$-solution, and functions as a local Euler-type approximation of the original system. Provided that the overall error (in time) has good concatenating properties, one can proceed with constructing a global quasi-constrained $\varepsilon$-solution (see Definition \ref{DefEpsSol} and Theorem \ref{ThAppSol}). This is a straight-forward machinery in the semi-groupal (or mild) formulation. For porous media equations, however, one works with $\Hmo{i}$-valued solutions which makes the estimates (in an Euler scheme, for instance) particularly sensitive to the $\Lo{2}{i}$-regularity of initial data and, globally, to time-integrated $\Lo{2}{i}$-norms. To understand the differences, the reader is invited to take a look at the constants appearing in Proposition \ref{Prop1}, inequalities \eqref{E3} and \eqref{E3'} and the induced impact in Proposition \ref{Prop2} and the necessary condition in Theorem \ref{ThNec}. As a consequence, the estimates on conditional expectations in Theorem \ref{ThNec} are, in general, not as good as the optimal speed ($t$ for deterministic classical systems leading to a choice of $\lambda=0$). These aspects motivate the further analysis in Subsection \ref{SubsecNec4} for finite-dimensional projections.

Finally, let us point out that for finite-dimensional systems, in the spirit of set-valued arguments, viability methods are often employed to provide epigraph-type comparisons, e.g. for viscosity solutions in \cite{FrankowskaPlaskacz2000} or, simply for forward or backward stochastic solutions, e.g. in \cite{Shi_2021}.  Of course, this kind of application relies on the presence of several components to be compared and it provides a further motivation for our partially coupled system.}

Bearing this in mind, at the end of the paper, we provide an application concerning the long-time behaviour of the solution of the stochastic porous media equation (see Section \ref{SecAppl} below). More precisely, considering a special case for $\beta_2,f_2$ in the equation \eqref{Eq00}, and applying the viability result obtained in this paper, for a particular set, we deduce an asymptotic exponential stabilization-type result, in projection, for the stochastic media equation. This can be put in connection with the controllability in projection theory  in \cite{Shyr}. However, we emphasize that, in contrast with \cite{Shyr},  here, the leading operator is fully-nonlinear. This makes the  problem very demanding, since the classical technique of linearization together with the tools from the theory of stabilization associated to linear systems (see the monograph \cite{Barbu_Sp}) cannot be used in the present context. To the best of our knowledge, this work provides for the first time a stabilization-type result for the stochastic porous media equation. From the practical point of view, this means that we fix a targeted moisture level, then design an internal controller which assures that, exponentially fast, the solution of the porous media equation approaches the targeted moisture value. Clearly, besides the theoretical value, this type of result is of high interest in stopping the cliff erosion, since once the moisture is confined to a certain values interval, one may prevent its erosion and degradation.\\

{The paper is organized as follows. Section \ref{SecStatement} presents the main assumptions used throughout the paper, the Euler-type scheme referred to as \textit{fundamental solution} in \eqref{EqFund} and the main notions of \emph{near-viability} and \emph{$\lambda$-quasi-tangency}. We wish to emphasize that, due to the specificities of porous media equations, distinction is made between $\mathbb{L}^2\pr{\mathcal{O}_i}$ and $\Hmo{i}$-notions.  The need for such distinction is clear from Section \ref{SecNec} (see also the aforementioned details on the originality of our work). Subsection \ref{SubsecNec1} gathers the integrability properties of the fundamental solution \eqref{EqFund} and the sharp estimates of right-continuity (see the constants in \eqref{E3} and \eqref{E3'}). Moreover, Proposition \ref{PropContInitialData} gives the dependence of initial data and perturbations for both solutions and fundamental solutions.  These tools concur in Subsection \ref{SubsecNec2} to reasonable estimates of the distance between the solution and its approximation via fundamental solutions. As a consequence, $\lambda$-quasi-tangency is shown to be necessary for viability in Theorem \ref{ThNec} (Subsection \ref{SubsecNec3}), for every $\lambda>0$. Subsection \ref{SubsecNec4} aims at sharpening $lambda$ to $0$ for finite-dimensional projections. Section \ref{SecSuff} deals with the sufficiency of the $0$-quasi-tangency condition in order to guarantee near viability. The sufficiency result is provided in Theorem \ref{ThSuf} by first constructing global approximating solution in Theorem \ref{ThAppSol}. Section \ref{SecAppl} contains the application of the theoretical results to obtain convenient conditions under which the porous media systems can be asymptotically driven towards the origin at an exponential speed.  Finally, the appendix presents, for the sake of completeness, the classical techniques of proof allowing to obtain the technical inequalities on solutions, fundamental solutions and approximating solutions.}\\

\textbf{Notations}\\
Throughout the paper, we will be using the following notations.
\begin{enumerate}
\item For an integer $n\in\mathbb{N}^*$, we let $\mathbb{R}^n$ stand for the standard $n$-dimensional Euclidean space. Euclidean norms will be denoted by $\abs{\cdot}$.
\item For a Lipschitz continuous function $\phi:\mathbb{R}^n\rightarrow\mathbb{R}^m$, we let $\pp{\phi}_1:=\underset{x\neq y}{\sup}\frac{\abs{\phi(x)-\phi(y)}}{\abs{x-y}}$ be its Lipschitz constant.
\item For $i\in\set{1,2}$ and $d_i\in\mathbb{N}^*$, we let $\mathcal{O}_i$ designate a bounded open subset of $\mathbb{R}^{d_i}$ with sufficiently regular ($C^2$) boundary $\partial\mathcal{O}_i$.
\begin{enumerate}
\item $\Lo{p}{i}$, $p\geq 1$ is the standard Banach space of real-valued $p$-power (Lebesgue-) integrable functions on $\mathcal{O}_i$.
\item $\Ho{i}$ is the space of $\mathbb{L}^2\pr{\mathcal{O}_i}$ functions that vanish on $\partial \mathcal{O}_i$ and such that the distributional derivative of first order belongs to $\mathbb{L}^p\pr{\mathcal{O}_i}$. The norm $\normi{\cdot}{1}{i}$ is given by \[\normi{\phi}{1}{i}^2:=\int_{\mathbb{R}^d}\abs{\nabla\phi\pr{\zeta}}^2d\zeta.\]
\item The dual of the aforementioned space is denoted by $H^{-1}\pr{\mathcal{O}_i}$. The associated norm is $\normi{\cdot}{-1}{i}
$ and the associated product $\scal{\cdot,\cdot}_{\Hmo{i}}$.
\item The dual of $\Lo{2}{i}$ is designated by $\Lo{2}{i}^*$ and the duality product by $\scal{\cdot,\cdot}_{\Lo{2}{i}^*,\Lo{2}{i}}$.
\end{enumerate}
{ 
\item We consider, for a generic Hilbert space $G$ a Wiener process $W$ on $G$. 
\item We let $\mathcal{L}_2$ designate the space of Hilbert-Schmidt operators on the spaces specified as arguments. We will often write $\mathcal{L}_2\pr{\Hmo{i}}$ instead of $\mathcal{L}_2\pr{G;\Hmo{i}}$.
\item Throughout the paper, $U$ is going to designate a compact (sub)space of a metric space (for simplicity, let us say a euclidean one $\mathbb{R}^d$. 
\item The stochastic problem is set on a complete filtered probability space $\pr{\Omega,\mathcal{F},\mathbb{F},\mathbb{P}}$ satisfying the usual assumptions of right-continuity and $\mathbb{P}$-completeness. 
\item The family of $U$-valued progressively measurable controls $u$ is going to be denoted by $\mathcal{U}$ and referred to as admissible control processes.
\item For random variables, integrability and measurability will be made explicit, e.g. , to designate $\Hmo{i}$-valued, $\mathcal{F}_t$-measurable, $\mathbb{P}$-square integrable random variables, we will use the notation $\mathbb{L}^2\pr{\Omega,\mathcal{F}_t,\mathbb{P};\Hmo{i}}$ for some $t>0$.}
\end{enumerate}
\section{Statement of the problem}\label{SecStatement}
{Throughout the paper and unless explicitly mentioned otherwise, the following standard assumption will hold true for $i\in\set{1,2}$.}
\begin{align}\label{Ass1}\begin{cases}
(i)\ &\beta_i:\mathbb{R}\rightarrow\mathbb{R}\textnormal{ are Lipschitz continuous and null at }0;\\[4pt]
(ii)\ &\textnormal{There exist } \alpha_i\geq 0\textnormal{ s.t. }\pr{\beta_i(r)-\beta_i(s)}\pr{r-s}\geq \alpha_i (r-s)^2,\ \forall r,s\in\mathbb{R};\\[4pt]
(iii)\ &f_i:\Hmo{1}\times \Hmo{2}\times U\rightarrow \Hmo{i},\\[4pt] &\forall x,x'\in\Hmo{1},\ y,y'\in \Hmo{2},\ u\in U,\\[4pt] 
&\norm{f_i(x,y,u)-f_i(x',y',u)}_{\Hmo{i}}\leq \pp{f_i}_{1}\pr{\norm{x-x'}_{\Hmo{1}}+\norm{y-y'}_{\Hmo{2}}};
\\[4pt]
&\sigma_i:\Hmo{1}\times \Hmo{2}\times U\rightarrow \mathcal{L}_2\pr{\Hmo{i}},\\[4pt]
 &\forall x,x'\in\Hmo{1},\ y,y'\in \Hmo{2},\ u\in U,\\[4pt] 
&\norm{\sigma_i(x,y,u)-\sigma_i(x',y',u)}_{\mathcal{L}_2\pr{\Hmo{i}}}\leq \pp{\sigma_i}_{1}\pr{\norm{x-x'}_{\Hmo{1}}+\norm{y-y'}_{\Hmo{2}}};\\[4pt]
(iv)\ &\underset{u\in U}{\sup}\norm{f_i(0,0,u)}_{\Hmo{i}}=:\norm{f_i}_0<\infty;\ \underset{u\in U}{\sup}\norm{\sigma_i(0,0,u)}_{\mathcal{L}_2\pr{\Hmo{i}}}=:\norm{\sigma_i}_0<\infty.
\end{cases}
\end{align}
\begin{remark}\label{Rem1}\begin{enumerate}
\item {We wish to emphasize that $\mathcal{L}_2\pr{\Hmo{i}}$ are used by making a slight abuse of notation. Since the equations have a common noise on some Hilbert space $G$, one should read $\mathcal{L}_2\pr{\Hmo{i}}:=\mathcal{L}_2\pr{G;\Hmo{i}}$.}
\item Note that, for a generic $\beta\in\set{\beta_1,\beta_2}$, this implies \begin{align*}\beta^2(r)\leq \pp{\beta}_1\beta(r)r.\end{align*}  Indeed, the previous inequality (ii) yields $\beta(r)\geq 0$, if $r\geq 0$.  As a consequence,  for $r\geq 0$, $\beta(r)=\abs{\beta(r)-\beta(0)}\leq \pp{\beta}_1r$. For $r\leq 0$, one reasons in the same way.
\item This can be easily extended to 
\begin{align}\label{Estim0}\overline \alpha \pr{\beta(r)-\beta(s)}^2\leq \pr{\beta(r)-\beta(s)}\pr{r-s}.\end{align}where $\overline \alpha:=\pr{\pp{\beta}_1+1}^{-1}>0$.
\item The assumptions also imply that, provided that $\Delta\beta_i(x)\in\Hmo{i}$, one has $\norm{\Delta\beta_i(x)}_{\Lo{2}{i}^*}\leq C\norm{x}_{\Lo{2}{i}}$. 
\end{enumerate}
\end{remark}
{ In addition, we impose the $\mathbb{L}^2$-regularity of the coefficients i.e.
\begin{align}\label{Ass2}\begin{cases}
(iii')&\forall x,x'\in\Lo{2}{1},\ y,y'\in \Lo{2}{2},\ u\in U,\\[4pt] 
&\norm{f_i(x,y,u)-f_i(x',y',u)}_{\Lo{2}{i}}\leq \pp{f_i}_{1}\pr{\norm{x-x'}_{\Lo{2}{1}}+\norm{y-y'}_{\Lo{2}{2}}};
\\[4pt]
&\norm{\sigma_i(x,y,u)-\sigma_i(x',y',u)}_{\mathcal{L}_2\pr{G;\Lo{2}{i}}}\leq \pp{\sigma_i}_{1}\pr{\norm{x-x'}_{\Lo{2}{1}}+\norm{y-y'}_{\Lo{2}{2}}};\\[4pt]
(iv')\ &\underset{u\in U}{\sup}\norm{f_i(0,0,u)}_{\Lo{2}{i}}<\infty;\ \underset{u\in U}{\sup}\norm{\sigma_i(0,0,u)}_{\mathcal{L}_2\pr{G;\Lo{2}{i}}}<\infty.
\end{cases}
\end{align}
}
{For initial times $t\geq 0$, we consider the following partially \footnote{"partially" coupled due to independent $\Delta \beta$ governing the components, see \cite{Shi_2021} for possible generalizations.}coupled porous media-type system. }
\begin{align}\label{Eq0}
\begin{cases}
dX^{t,\xi,\eta,u}(s)&=\pr{\Delta\beta_1\pr{X^{t,\xi,\eta,u}(s)}+f_1\pr{X^{t,\xi,\eta,u}(s),Y^{t,\xi,\eta,u}(s),u(s)}}ds\\[4pt]
&+\sigma_1\pr{X^{t,\xi,\eta,u}(s),Y^{t,\xi,\eta,u}(s),u(s)}dW(s),\\[4pt]
dY^{t,\xi,\eta,u}(s)&=\pr{\Delta\beta_2\pr{Y^{t,\xi,\eta,u}(s)}+f_2\pr{X^{t,\xi,\eta,u}(s),Y^{t,\xi,\eta,u}(s),u(s)}}ds\\[4pt] &+\sigma_2\pr{X^{t,\xi,\eta,u}(s),Y^{t,\xi,\eta,u}(s),u(s)}dW(s),\ s\geq t;\\[4pt]
X^{t,\xi,\eta,u}(t)&=\xi\in\mathbb{L}^2\pr{\Omega,\mathcal{F}_t,\mathbb{P};\Hmo{1}},\ Y^{t,\xi,\eta,u}(t)=\eta\in\mathbb{L}^2\pr{\Omega,\mathcal{F}_t,\mathbb{P};\Hmo{2}},
\end{cases}
\end{align}
together with the fundamental solution
\begin{align}\label{EqFund}
\begin{cases}
d\mathbb{X}^{t,\xi,\eta,u}(s)&=\pr{\Delta\beta_1\pr{\XX(s)}+f_1\pr{\xi,\eta,u(s)}}ds+\sigma_1\pr{\xi,\eta,u(s)}dW(s),\\[4pt]
d\mathbb{Y}^{t,\xi,\eta,u}(s)&=\pr{\Delta\beta_2\pr{\YY(s)}+f_2\pr{\xi,\eta,u(s)}}ds+\sigma_2\pr{\xi,\eta,u(s)}dW(s),\ s\geq t;\\[4pt]
\XX(t)&=\xi\in\mathbb{L}^2\pr{\Omega,\mathcal{F}_t,\mathbb{P};\Hmo{1}}, \ \YY(t)=\eta\in\mathbb{L}^2\pr{\Omega,\mathcal{F}_t,\mathbb{P};\Hmo{2}}.
\end{cases}
\end{align}
These equations are controlled with a progressively measurable process $u$ taking its values in a compact metric space $U$. The existence and uniqueness of solutions for the (coupled) system \eqref{Eq0} (resp. \eqref{EqFund}) follows in a standard way (see \cite[Chapter 2]{BDPR_2016}, or the recent \cite{HLL2021} for more general coupled systems).
In connection with these equations, we consider the following concepts.
\begin{definition}\label{DefNearViab}
Given a closed set $K\subset \Hmo{1}\times \Hmo{2}$,  it is said to be \begin{enumerate}
\item \textbf{$\mathbb{L}^2$-nearly viable} with respect to \eqref{Eq0} on the finite interval $\pp{0,T>0}$ if, for every initial time $t\in\left[0,T\right)$ and every initial pair $\pr{\xi,\eta}\in \mathbb{K}_t:=\mathbb{L}^2\pr{\Omega,\mathcal{F}_t,\mathbb{P};K}$, such that $\pr{\xi,\eta}\in \Lo{2}{1}\times\Lo{2}{2}$, $\mathbb{P}$-a.s., we have \begin{equation}
\label{EqNearViab}
\inf_{u\in \mathcal{U}}\sup_{s\in\pp{t,T}}d\pr{\pr{\X(s),\Y(s)},\mathbb{K}_s{ \cap\mathbb{L}^2\pr{\Omega,\mathcal{F}_s,\mathbb{P};\Lo{2}{1}\times\Lo{2}{2}}}}=0.
\end{equation}{The distance $d$ is meant in $\mathbb{L}^2\pr{\Omega,\mathcal{F},\mathbb{P};\Hmo{1}\times\Hmo{2}}$.}
\item \textbf{nearly viable} with respect to \eqref{Eq0} on the finite interval $\pp{0,T>0}$ if, for every initial time $t\in\left[0,T\right)$ and every initial pair $\pr{\xi,\eta}\in \mathbb{K}_t:=\mathbb{L}^2\pr{\Omega,\mathcal{F}_t,\mathbb{P};K}$, we have \begin{equation}
\label{EqNearViab'}
\inf_{u\in \mathcal{U}}\sup_{s\in\pp{t,T}}d\pr{\pr{\X(s),\Y(s)},\mathbb{K}_s}=0.
\end{equation}
\end{enumerate}
\end{definition}
\begin{remark}
\begin{enumerate}
\item If an optimal control in \eqref{EqNearViab'} exists, then near viability is equivalent with viability i.e.  $\pr{\X(s),\Y(s)}\in K$,  $\mathbb{P}$-a.s. and for all $s\in\pp{t,T}$. 
\item {The distance in \eqref{EqNearViab'} can be replaced with the distance in $\Hmo{1}\times\Hmo{2}$ from $\X(s)$ to $K$ whenever the orthogonal projector set-valued function is non-empty (cf. \cite[Corollary 8.2.13]{AubF90}).}
\end{enumerate}
\end{remark}
{The second concept is a \emph{quasi-tangency} one. It describes the adequacy of the fundamental solution as a candidate to state-constrained Euler schemes. In particular, this notion quantifies the local behaviour of the distance to the set of constraints.}
\begin{definition}\label{DefLambdaQT}
Let us fix the finite interval $\pp{0,T>0}$ and $t\in\left[0,T\right)$.  A closed set $K\subset \Hmo{1}\times \Hmo{2}$ satisfies the \textbf{$\lambda$-quasi-tangency condition} with respect to the control system \eqref{Eq0} at $\pr{\xi,\eta}\in \mathbb{K}_t$ if 
\begin{equation*}
\begin{split} \underset{\varepsilon\rightarrow0+}{\lim\inf}\inf &\left\{ \frac{1}{\varepsilon}\mathbb{E}\pp{\norm{\XX(t+\varepsilon)-\theta_1}_{\Hmo{1}}^2+\norm{\YY(t+\varepsilon)-\theta_2}_{\Hmo{2}}^2} \right. \\[4pt] &+\frac{1}{\varepsilon^{2-2\lambda}}\mathbb{E}\pp{\norm{\ef{\XX(t+\varepsilon)-\theta_1}{t}}_{\Hmo{1}}^2}\\[4pt] &+\frac{1}{\varepsilon^{2-2\lambda}}\mathbb{E}\pp{\norm{\ef{\YY(t+\varepsilon)-\theta_2}{t}}_{\Hmo{2}}^2}:\\[4pt] & \left. u,\ U\textnormal{-valued,  progressively measurable},\ \theta=\pr{\theta_1,\theta_2}\in\mathbb{K}_t:=\mathbb{L}^2\pr{\Omega,\mathcal{F}_t,\mathbb{P};K} \right\} =0. 
\end{split}
\end{equation*} 

If this condition holds true for every $t\in\left[0,T\right)$ and every $\pr{\xi,\eta}\in \Hmo{1}\times \Hmo{2}$, then we will simply say that $K$ satisfies the $\lambda$-quasi-tangency condition.\\
{ When $\pr{\xi,\eta}\in \mathbb{K}_t$ have $\Lo{2}{1}\times\Lo{2}{2}$ regularity, we will require the same regularity for $\pr{\theta_1,\theta_2}$.  In some sense, one can talk about the \textbf{$\lambda$-quasi-tangency condition in $\mathbb{L}^2$}. For further connections, please refer to Theorem \ref{ThNec}.}
\end{definition}
\section{The necessary condition}\label{SecNec}
{From now on, we let $C>0$ be a generic constant only depending on the fixed time horizon $T>0$ and the dynamical characteristics specified in the assumptions but independent of initial conditions and changing/asymptotic parameters i.e. $\varepsilon\rightarrow 0, \ n\rightarrow\infty$ and so on. The constant is allowed to change from one line to another in order to simplify the presentation.}
\subsection{Basic regularity and initial estimates}\label{SubsecNec1}
{We gather here some useful regularity properties for the fundamental solution $\XX$ and the dependency of $\X$ on initial data.}
\begin{proposition}\label{Prop1}Let the time horizon $T>0$ be fixed and Assumptions \ref{Ass1} and \ref{Ass2} to hold true. Then, there exists a constant $C>0$ (generic, dependent of $T$ but not of the initial data) such that the following hold true.
\begin{enumerate}
\item For every $\xi\in\mathbb{L}^2\pr{\Omega,\mathcal{F}_t,\mathbb{P};\Hmo{1}}$ and every $\eta\in\mathbb{L}^2\pr{\Omega,\mathcal{F}_t,\mathbb{P};\Hmo{2}}$, and every $t\leq s\leq T$,
\begin{equation}\label{L2estim}\begin{split}
&\mathbb{E}\pp{\sup_{t\leq r\leq s}\norm{\XX(r)}_{\Hmo{1}}^2+\sup_{t\leq r\leq s}\norm{\YY(r)}_{\Hmo{2}}^2}\\[4pt] &+\mathbb{E}\pp{\int_t^s\pr{\norm{\beta_1\pr{\XX(l)}}_{\Lo{2}{1}}^2+\norm{\beta_2\pr{\YY(l)}}_{\Lo{2}{2}}^2}dl}
\\[4pt] &\leq C\pr{1+\mathbb{E}\pp{\norm{\xi}_{\Hmo{1}}^2+\norm{\eta}_{\Hmo{2}}^2}}.
\end{split}\end{equation}Furthermore, let us assume that\begin{equation}
\label{Ass+}
\textnormal{for each }i\in\set{1,2},\textnormal{ either }\beta_i\equiv 0\textnormal{ or }\alpha_i>0. 
\end{equation}Then we have the following.
\item For every $\xi\in\mathbb{L}^2\pr{\Omega,\mathcal{F}_t,\mathbb{P};\Lo{2}{1}}$ and every $\eta\in\mathbb{L}^2\pr{\Omega,\mathcal{F}_t,\mathbb{P};\Hmo{2}}$, and every $t\leq s\leq T$,\begin{equation}\label{E3}\begin{split}
\mathbb{E}\pp{\sup_{t\leq r\leq s}\norm{\XX(r)-\xi}_{\Hmo{1}}^2}\leq C\pr{1+\mathbb{E}\pp{\norm{\xi}_{\Lo{2}{1}}^2+\norm{\eta}_{\Hmo{2}}^2}}(s-t).
\end{split}
\end{equation}
\item For every $\xi\in\mathbb{L}^2\pr{\Omega,\mathcal{F}_t,\mathbb{P};\Hmo{1}}$ and every $\eta\in\mathbb{L}^2\pr{\Omega,\mathcal{F}_t,\mathbb{P};\Lo{2}{2}}$, and every $t\leq s\leq T$,\begin{equation}\label{E3'}\begin{split}
\mathbb{E}\pp{\sup_{t\leq r\leq s}\norm{\YY(r)-\eta}_{\Hmo{2}}^2}\leq C\pr{1+\mathbb{E}\pp{\norm{\eta}_{\Lo{2}{2}}^2+\norm{\xi}_{\Hmo{1}}^2}}(s-t).
\end{split}
\end{equation}
\end{enumerate}
\end{proposition}
{The proof relies on standard considerations: Itô's formula for square norms $\normi{\cdot}{-1}{i}^2$, the monotonicity of $\beta_i$, usual Burkholder-Davis-Gundy inequalities, and Gronwall's lemma. For our readers' sake, the proofs of these natural estimates will be relegated to Section \ref{SecAppendix}.}
{
\begin{remark}\label{RemSec}A careful look at the estimates in the proof (in \eqref{E1}) shows that we can actually get a sharper bound involving \begin{equation}\label{estimExp}
\mathbb{E}\pp{\int_t^s\norm{\XX-\xi}_{\Lo{2}{1}}^2dl}\leq C\pr{1+\mathbb{E}\pp{\norm{\xi}_{\Lo{2}{1}}^2+\norm{\eta}_{\Hmo{2}}^2}}(s-t).
\end{equation}
This bound is superfluous for the theoretical proofs, but might be interesting in applications (see Section \ref{SecAppl}.
\end{remark}
}
At this point, let us emphasize the continuous dependence of the initial data.
\begin{proposition}\label{PropContInitialData}Let the time horizon $T>0$ be fixed and Assumptions \eqref{Ass1} {and \eqref{Ass2}} to hold true. Then, there exists a constant $C>0$ such that, for every $\xi,  \xi'\in\Hmo{1}$ and every $\eta,\eta'\in\Hmo{2}$ and every $t\leq s\leq T$,\begin{equation*}
\begin{split}
&\mathbb{E}\pp{\sup_{t\leq r\leq s}\normi{\X(r)-\Xp(r)}{-1}{1}^2+\int_t^s\norm{\beta_1\pr{\X(l)}-\beta_1\pr{\Xp(l)}}_{\Lo{2}{1}}^2dl}\\[4pt]
&+\mathbb{E}\pp{\sup_{t\leq r\leq s}\normi{\Y(r)-\Yp(r)}{-1}{2}^2+\int_t^s\norm{\beta_2\pr{\Y(l)}-\beta_2\pr{\Yp(l)}}_{\Lo{2}{2}}^2dl}\\[4pt]&\leq C\pr{\norm{\xi-\xi'}_{\Hmo{1}}^2+\norm{\eta-\eta'}_{\Hmo{2}}^2},
\end{split}\end{equation*} and\begin{equation*}
\begin{split}
&\mathbb{E}\pp{\sup_{t\leq r\leq s}\normi{\XX(r)-\XXp(r)}{-1}{1}^2+\int_t^s\norm{\beta_1\pr{\XX(l)}-\beta_1\pr{\XXp(l)}}_{\Lo{2}{1}}^2dl}\\[4pt]
&+\mathbb{E}\pp{\sup_{t\leq r\leq s}\normi{\YY(r)-\YYp(r)}{-1}{2}^2+\int_t^s\norm{\beta_2\pr{\YY(l)}-\beta_2\pr{\YYp(l)}}_{\Lo{2}{2}}^2dl}\\[4pt]&\leq C\pr{\norm{\xi-\xi'}_{\Hmo{1}}^2+\norm{\eta-\eta'}_{\Hmo{2}}^2}.
\end{split}\end{equation*}
\end{proposition}
{The proof is, again, based on standard techniques and will be relegated to Section \ref{SecAppendix} for our readers' comfort.}
\subsection{Estimates in average in $\Hmo{\cdot}$. The necessity}\label{SubsecNec2}
We now strive to show that $\pr{\XX,\YY}$ provides a good replacement for the actual solution $\pr{\X,\Y}$. {Again, the proof is postponed to the Appendix and it relies on classical techniques Itô's formula, monotonicity of the driving$\beta_i$, the Lispchitz continuity for the coefficients and Burkolder-Davis-Gundy and Gronwall inequalities. }
\begin{proposition}\label{Prop2}Let us assume \eqref{Ass1} and \eqref{Ass+} to hold true. Then, there exists a constant $C>0$ (generic, able to change from one line to another and only depending on the time horizon $T>0$ and the coefficients) such that, for every $\xi\in\mathbb{L}^2\pr{\Omega,\mathcal{F}_t,\mathbb{P};\Lo{2}{1}}$, every $\eta\in\mathbb{L}^2\pr{\Omega,\mathcal{F}_t,\mathbb{P};\Lo{2}{2}}$,  every admissible control $u\in\mathcal{U}$, and every $t\leq s\leq T$,
\begin{equation}
\label{FirstEstim}\begin{split}
&\mathbb{E}\pp{\sup_{t\leq r\leq s}\norm{\XX(r)-\X(r)}_{\Hmo{1}}^2+\int_t^s\norm{\XX(l)-\X(l)}_{\Lo{2}{1}}dl}\\[4pt]
&+\mathbb{E}\pp{\sup_{t\leq r\leq s}\norm{\YY(r)-\Y(r)}_{\Hmo{2}}^2+\int_t^s\norm{\YY(l)-\Y(l)}_{\Lo{2}{2}}dl}\\[4pt]
&\leq C\pr{1+\mathbb{E}\pp{\norm{\xi}_{\Lo{2}{1}}^2+\norm{\eta}_{\Lo{2}{2}}^2}}(s-t)^{2}.
\end{split}
\end{equation}
\end{proposition}
\subsection{The necessary condition}\label{SubsecNec3}
These estimates allow us to prove the following necessity condition.
\begin{theorem}\label{ThNec}
Let $K\subset \Hmo{1}\times\Hmo{2}$ be a closed set. 
\begin{enumerate}
\item If $K$ is $\mathbb{L}^2$-nearly viable with respect to \eqref{Eq0}, then $K$ enjoys the $\lambda$-quasi tangency condition for every $\lambda > 0$ {(in $\mathbb{L}^2$, cf. Definition \ref{DefLambdaQT})} with respect to \eqref{Eq0} at every $\pr{\xi,\eta}\in\Lo{2}{1}\times\Lo{2}{2}$.
\item Furthermore,  let us assume that $\mathbb{L}^2\pr{\Omega,\mathcal{F}_t,\mathbb{P};\Lo{2}{1}\times\Lo{2}{2}}$ is relatively dense in $\mathbb{K}_t$ for some $T>0$ and each $t\in\left[0,T\right)$. Then \begin{enumerate}
\item $K$ is nearly viable with respect to \eqref{Eq0} if and only if it is $\mathbb{L}^2$-nearly viable with respect to \eqref{Eq0} (on $\pp{0,T}$);
\item $K$ satisfies the $\lambda$-quasi tangency condition with respect to \eqref{Eq0} if and only if it satisfies this condition for couples $\pr{\xi,\eta}\in \mathbb{L}^2\pr{\Omega,\mathcal{F}_t,\mathbb{P};\Lo{2}{1}\times\Lo{2}{2}}\cap\mathbb{K}_t$.
\end{enumerate}
\end{enumerate}
\end{theorem}
\begin{remark}
Before proceeding with the proof, we point out that the density requirement is quite weak. Indeed, $\Lo{2}{i}$ is continuously and densely embedded into $\Hmo{i}$, such that this requirement comes down to asking $K$ to be the closure of its interior, and some metric considerations on $\Omega$ (Polish space as it is usually the case with the canonical Wiener constructions) in order to be able to approximate measurable functions with "piecewise" constant ones.
\end{remark}
\begin{proof}
Let us fix the time horizon $T>0$ and the initial time $t\in\left[0,T\right)$.\\
1. For the first assertion, we also fix $\varepsilon>0$ and $\pr{\xi,\eta}\in \mathbb{K}_t{ \cap\mathbb{L}^2\pr{\Omega,\mathcal{F}_t,\mathbb{P};\Lo{2}{1}\times\Lo{2}{2}}}$. Then, there exists an admissible control $u\in\mathcal{U}$ (progressively measurable, $U$-valued) such that \[d^2\pr{\pr{\X(s),\Y(s)},\mathbb{K}_s{\cap\mathbb{L}^2\pr{\Omega,\mathcal{F}_s,\mathbb{P};\Lo{2}{1}\times\Lo{2}{2}}}}<\varepsilon^3,\ \forall s\in\pp{t,T}.\]One is able to pick some $\theta=\pr{\theta_1,\theta_2}\in\mathbb{K}_{t+\varepsilon}{ \cap\mathbb{L}^2\pr{\Omega,\mathcal{F}_{t+\varepsilon},\mathbb{P};\Lo{2}{1}\times\Lo{2}{2}}}$ such that \[\mathbb{E}\pp{\normi{\X(t+\varepsilon)-\theta_1}{-1}{1}^2+\normi{\Y(t+\varepsilon)-\theta_2}{-1}{2}^2}\leq \varepsilon^3.\]Owing to Proposition \ref{Prop2}, it follows from the later inequality that \begin{align*}
&\mathbb{E}\pp{\normi{\XX(t+\varepsilon)-\theta_1}{-1}{1}^2+\normi{\YY(t+\varepsilon)-\theta_1}{-1}{1}^2}\\[4pt]&\leq C\pr{1+\mathbb{E}\pp{\norm{\xi}_{\Lo{2}{1}}^2+\norm{\eta}_{\Lo{2}{2}}^2}}\varepsilon^2.
\end{align*}The conclusion follows.\\
2. Due to the density assumption, for every $\pr{\xi',\eta'}\in \mathbb{K}_t$ and $\varepsilon>0$, one is able to find $\pr{\xi,\eta}\in\mathbb{K}_t\cap\mathbb{L}^2\pr{\Omega,\mathcal{F}_t,\mathbb{P};\Lo{2}{1}\times\Lo{2}{2}}$ such that \[\mathbb{E}\pp{\norm{\xi'-\xi}_{\Hmo{1}}^2+\norm{\eta'-\eta}_{\Hmo{2}}^2}\leq \varepsilon^3.\]
(a) If $K$ is $\mathbb{L}^2$-nearly viable, then, there exists $u\in\mathcal{U}$ such that \[d_{\mathbb{L}^2\pr{\Omega,\mathcal{F},\mathbb{P};\Hmo{1}\times\Hmo{2}}}\pr{\pr{\X(s),\Y(s)},\mathbb{K}_s}\leq \frac{1}{2}\varepsilon,\ \forall s\in\pp{t,T}.\]Then, owing to Proposition \ref{PropContInitialData}, it follows that, for all $s\in \pp{t,T}$, \begin{align*}
&d_{\mathbb{L}^2\pr{\Omega,\mathcal{F},\mathbb{P};\Hmo{1}\times\Hmo{2}}}\pr{\pr{\Xp(s),\Yp(s)},\mathbb{K}_s}\\
\leq &d_{\mathbb{L}^2\pr{\Omega,\mathcal{F},\mathbb{P};\Hmo{1}\times\Hmo{2}}}\pr{\pr{\X(s),\Y(s)},\mathbb{K}_s}\\&+\sqrt{\normi{\X(s)-\Xp(s)}{-1}{1}^2+\normi{\Y(s)-\Yp(s)}{-1}{2}^2}\leq \varepsilon,
\end{align*}for $\varepsilon$ small enough chosen to compensate the generic constant in Proposition \ref{PropContInitialData}. It follows that $K$ is nearly viable w.r.t. \eqref{Eq0}.\\
(b) The quasi-tangency is proven in a similar way from the quantities involved in Definition \ref{DefLambdaQT} and using the estimates on $\normi{\XX-\XXp}{-1}{1}^2+\normi{\YY-\YYp}{-1}{2}^2$ from Proposition~ \ref{PropContInitialData}.
\end{proof}
\subsection{Estimating the conditional expectations in $\Hmo{\cdot}$}\label{SubsecNec4}
While the previous result provides us with a good speed for the stochastic case,  one would expect a better estimate for the "deterministic" part.  In particular, one is interested in cases when (much like the finite-dimensional setting), the quasi-tangency can be obtained with $\lambda=0$.\\

To this purpose, we look into $\normi{\ef{\XX-\X}{t}}{-1}{1}^2$ and the analogous part for the second component. {We wish to point out that a direct (iterated) method works for deterministic equations. A contribution in this setting \cite{CGM_2022+} will be available soon, but it fails to provide a reasonable estimate for the stochastic case.}\\
For $i\in\set{1,2}$, let us denote by \begin{enumerate}
\item $\lambda^i_j$ the eigenvalues of $-\Delta$ w.r.t. $\mathcal{O}_i$ and ordered increasingly;
\item $\tilde e_j^i$ the $\Lo{2}{i}$ orthonormal basis of eigen-functions vanishing on $\partial\mathcal{O}_i$ and corresponding to these eigenvalues. It is clear that $e_j^i:={\sqrt{\lambda_j^i}}\tilde e_j^i$ is orthonormal in $\Hmo{i}$;
\item For $j\geq 1$, we let $\Pi^i_j:\Hmo{i}\rightarrow\Lo{2}{i}\subset\Hmo{i}$ be defined by \begin{equation}
\label{Pij}
\Pi^i_j(x)=\sum_{1\leq k\leq j}\scal{x,e_k^i}_{\Hmo{i}}e_k^i.
\end{equation}
The reader is invited to note this is an orthogonal projector in $\Hmo{i}$ and, in particular, $\norm{\Pi_i^j(x)}_{\Hmo{i}}\leq \norm{x}_{\Hmo{i}}$.
\end{enumerate}
In particular, \begin{enumerate}
\item If $x\in\Lo{2}{i}$, then \begin{equation*}
\scal{\Delta x,e_j^i}_{\Hmo{i}}={-\sqrt{\lambda_j^i}}\scal{x, \tilde e_j^1}_{\Lo{2}{i}}; \; \; \; \;
\scal{x,e_j^i}_{\Hmo{i}}=\frac{1}{\sqrt{\lambda_j^i}}\scal{x,\tilde e_j^1}_{\Lo{2}{i}};
\end{equation*}
\item As a consequence,  $\Pi^i_j(x)=\sum_{1\leq k\leq j}\scal{x,\tilde e_k^i}_{\Lo{2}{i}}\tilde e_k^i$ is the usual $\Lo{2}{i}$ projection onto\\ $span_{\Lo{2}{1}}\set{\tilde e_k^i:\ 1\leq k\leq j}$ when restricted to $\Lo{2}{i}$.
\item 
We have
\begin{equation*}
\begin{split}
\norm{\Pi_j^i\pr{\Delta x}}_{\Lo{2}{i}}^2 &=\underset{1\leq k\leq j}{\sum}\pr{\lambda_k^i}^2\scal{x, \tilde e_k^i}_{\Lo{2}{1}}^2\\[4pt] 
&\leq \underset{1\leq k\leq j}{\max}\pr{\lambda_k^i}^2\norm{\Pi_j^i(x)}_{\Lo{2}{1}}^2=\pr{\lambda_j^i}^2\norm{\Pi_j^i(x)}_{\Lo{2}{i}}^2;
\end{split}
\end{equation*}
\item 
Moreover
\begin{equation*}
\begin{split}\normi{\Pi_j^i\pr{\Delta x}}{-1}{i}^2&=\underset{1\leq k\leq j}{\sum}\pr{\lambda_k^i}^2\scal{x,e_k^i}_{\Hmo{i}}^2\\[4pt] &\leq \pr{\underset{1\leq k\leq j}{\max}\pr{\lambda_k^i}^2}\normi{\Pi_j^i(x)}{-1}{i}^2;
\end{split}
\end{equation*}
\item 
Furthermore
\begin{equation*}
\begin{split}
\norm{\Pi_j^i\pr{x}}_{\Lo{2}{i}}^2&=\underset{1\leq k\leq j}{\sum}\scal{x,\tilde e_k^i}_{\Lo{2}{i}}^2\\[4pt]
&=\underset{1\leq k\leq j}{\sum}{\lambda_k^i}\scal{x,e_k^i}_{\Hmo{i}}^2\leq \pr{\underset{1\leq k\leq j}{\max}{\lambda_k^i}}\norm{\Pi_j^i(x)}_{\Hmo{i}}^2;
\end{split}
\end{equation*}
\item 
In addition
\begin{equation*}
\begin{split}
\norm{\Pi_j^i(x)}_{\Hmo{i}}^2&=\underset{1\leq k\leq j}{\sum}\scal{x,e_k^i}_{\Hmo{i}}^2\\[4pt]&=\underset{1\leq k\leq j}{\sum}\frac{1}{\lambda_k^i}\scal{x,\tilde e_k^i}_{\Lo{2}{i}}^2\leq \frac{1}{\underset{1\leq k\leq j}{\min}{\lambda_k^i}}\norm{\Pi_j^i\pr{x}}_{\Lo{2}{i}}^2.
\end{split}
\end{equation*}
\item Note also that $\Pi_j^i$ commutes with $\Delta$.
\end{enumerate}
\begin{proposition}
Let us assume \eqref{Ass1} and \eqref{Ass+} to hold true and let us fix $j,j'\in\mathbb{N}^*$. Then, there exists a constant $C>0$ (generic, able to change from one line to another and only depending on the time horizon $T>0$ and the coefficients) such that, for every $\xi\in\mathbb{L}^2\pr{\Omega,\mathcal{F}_t,\mathbb{P};\Lo{2}{1}}$ and every $\eta\in\mathbb{L}^2\pr{\Omega,\mathcal{F}_t,\mathbb{P};\Lo{2}{2}}$, and every $t\leq s\leq T$,
\begin{equation}
\label{SecondEstim1}\begin{split}
&\mathbb{E}\pp{\sup_{t\leq r\leq s}\normi{\ef{\Pi_j^1\pr{\XX(r)-\X(r)}}{t}}{-1}{1}^2}\\[4pt]
&+\mathbb{E}\pp{\sup_{t\leq r\leq s}\normi{\ef{\Pi_{j'}^2\pr{\YY(r)-\Y(r)}}{t}}{-1}{2}^2}\\[4pt]
&\leq C\pr{1+\pr{\frac{\underset{1\leq k\leq j}{\max}\lambda_k^1}{\underset{1\leq k\leq j}{\min}\sqrt{\lambda_k^1}}}^2\vee\pr{\frac{\underset{1\leq k\leq j'}{\max}\lambda_k^2}{\underset{1\leq k\leq j'}{\min}\sqrt{\lambda_k^2}}}^2}\pr{1+\mathbb{E}\pp{\norm{\xi}_{\Lo{2}{1}}^2+\norm{\eta}_{\Lo{2}{2}}^2}}(s-t)^{3}.
\end{split}
\end{equation}
\end{proposition}
\begin{proof}
One easily notes that \begin{equation*}
\begin{split}
&\pr{\delta_j^1(r)}^2:=\normi{\ef{\Pi_j^1\pr{\XX(r)-\X(r)}}{t}}{-1}{1}^2\\[4pt]
&\begin{split}=&\left\langle \ef{\int_t^r\Pi^1_j\pr{\Delta\beta_1\pr{\XX(l)}-\Delta\beta_1\pr{\X(l)}}dl}{t}, \right. \\[4pt]& 
\left. \ef{\Pi_j^1\pr{\XX(r)-\X(r)}}{t}\right\rangle_{\Hmo{1}} \end{split}
\\[4pt]
&\begin{split}+&
\left\langle\ef{\int_t^r\Pi_j^1\pr{f_1\pr{\xi,\eta,u(l)}-f_1\pr{\X(l),\Y(l),u(l)}}dl}{t}\right.,
\\[4pt]&
\left.\ef{\Pi_j^1\pr{\XX(r)-\X(r)}}{t}\right\rangle_{\Hmo{1}}\end{split}\\[4pt]
\leq &\ef{\int_t^r\normi{\Pi^1_j\pr{\Delta\beta_1\pr{\XX(l)}-\Delta\beta_1\pr{\X(l)}}}{-1}{1}dl}{t}\delta_j^1(r)\\[4pt]&+\ef{\int_t^r\normi{\Pi^1_j\pr{f_1\pr{\xi,\eta,u(l)}-f_1\pr{\X(l),\Y(l),u(l)}}}{-1}{1}dl}{t}\delta_j^1(r)\\[4pt]
\leq &\frac{\underset{1\leq k\leq j}{\max}\lambda_k^1}{\underset{1\leq k\leq j}{\min}\sqrt{\lambda_k^1}}\pp{\beta}_1\ef{\int_t^r\norm{\XX(l)-\X(l)}_{\Lo{2}{1}}dl}{t}\delta_j^1(r)\\[4pt]&+\pp{f}_1(r-t)\ef{\sup_{t\leq l\leq r}\pr{\normi{\XX(l)-\xi}{-1}{1}+\normi{\YY-\eta}{-1}{2}}}{t}\delta_j^1(r)
\end{split}
\end{equation*}
As a consequence,  using classical inequalities $ab\leq \frac{1}{8}b^2+2a^2$ and Jensen's inequality, one gets
\begin{align*}\pr{\delta_j^1(r)}^2\leq &C\pr{1+\pr{\frac{\underset{1\leq k\leq j}{\max}\lambda_k^1}{\underset{1\leq k\leq j}{\min}\sqrt{\lambda_k^1}}}^2}\ef{\pr{\int_t^r\norm{\pr{\XX(l)-\X(l)}}_{\Lo{2}{1}}dl}^2}{t}\\[4pt]&+C(r-t)^2\ef{\pr{\sup_{t\leq l\leq r}\normi{\XX(l)-\xi}{-1}{1}^2+\sup_{t\leq l\leq r}\normi{\YY-\eta}{-1}{2}^2}}{t}\\[4pt]
\leq &C\pr{1+\pr{\frac{\underset{1\leq k\leq j}{\max}\lambda_k^1}{\underset{1\leq k\leq j}{\min}\sqrt{\lambda_k^1}}}^2}(r-t)\ef{\int_t^r\norm{\pr{\XX(l)-\X(l)}}_{\Lo{2}{1}}^2dl}{t}\\[4pt]&+C(r-t)^2\ef{\pr{\sup_{t\leq l\leq r}\normi{\XX(l)-\xi}{-1}{1}^2+\sup_{t\leq l\leq r}\normi{\YY-\eta}{-1}{2}^2}}{t}.
\end{align*}
The conclusion follows by taking $\underset{t\leq r\leq s}{\sup}$ and invoking Proposition \ref{Prop2} to deal with the first term and Proposition \ref{Prop1} for the remaining term(s).
\end{proof}

By putting $s=t+\varepsilon$ in the previous proposition, one gets the stronger quasi-tangency condition.
\begin{corollary}
Let us assume that there exist $j,j'\geq 1$ such that $K\subset \Pi^1_j\Hmo{1}\times \Pi^2_{j'}\Hmo{2}$. If $K$ satisfies the near-viability property with respect to \eqref{Eq0} on $\pp{0,T>0}$, then
\begin{enumerate}
\item it is quasi-tangent with $\lambda=0$,  for every initial $t\in\left[0,T\right)$ and every \\$\pr{\xi,\eta}\in\mathbb{L}^2\pr{\Omega,\mathcal{F}_t,\mathbb{P};\Lo{2}{1}\times \Lo{2}{2}}\cap\mathbb{K}_t$;
\item If, in addition, $\mathbb{L}^2\pr{\Omega,\mathcal{F}_t,\mathbb{P};\Lo{2}{1}\times\Lo{2}{2}}$ is relatively dense in $\mathbb{K}_t$ for every, $t\in\left[0,T\right)$, then $K$ is quasi-tangent with $\lambda=0$,  for every initial $t\in\left[0,T\right)$ and every $\pr{\xi,\eta}\in\mathbb{K}_t$.
\end{enumerate}
\end{corollary}
{
\section{Sufficiency}\label{SecSuff}
Let us note that, for $\lambda\in\pp{0,\frac{1}{2}}$, the $\lambda$-quasi-tangency condition can, alternatively, be expressed in a sequential formulation
\begin{proposition}\label{PropSeqQT}
The set $K$ has the $\lambda$-quasi-tangency property with respect to \eqref{Eq0} at time $t\in\left[0,T\right)$ at the point $\pr{\xi,\eta}\in \mathbb{K}_t$ if and only if the following hold simultaneously.\begin{enumerate}\item There exists a sequence $\pr{\varepsilon_n}_{n\geq 1}\subset\mathbb{R}_+$ such that $\underset{n\rightarrow\infty}{\lim}\varepsilon_n=0$ and there exist sequences $a^i_n,b^i_n\in \mathbb{L}^2\pr{\Omega,\mathcal{F}_{t+\varepsilon_n},\mathbb{P};\Hmo{i}}$, for $i\in \set{1,2}$ such that
\begin{enumerate}
\item $\underset{n\rightarrow\infty}{\lim}\mathbb{E}\pp{\normi{a_n^i}{-1}{i}^2}=\underset{n\rightarrow\infty}{\lim}\mathbb{E}\pp{\normi{b_n^i}{-1}{i}^2}=0$;
\item $\ef{b_n^i}{t}=0,\ \mathbb{P}$-a.s.  and $a_n^i$ is $\mathcal{F}_t$-measurable;
\end{enumerate}
and
\item there exists 
$u_n\in\mathcal{U}$ such that \begin{equation*}
\pr{\mathbb{X}^{t,\xi,\eta,u_n}\pr{t+\varepsilon_n}+\varepsilon_n^{1-\lambda}a_n^1+\sqrt{\varepsilon_n}b_n^1,\mathbb{Y}^{t,\xi,\eta,u_n}\pr{t+\varepsilon_n}+\varepsilon_n^{1-\lambda}a_n^2+\sqrt{\varepsilon_n}b_n^2}\in \mathbb{K}_{t+\varepsilon},\ \mathbb{P}\textnormal{-a.s.}
\end{equation*}
\end{enumerate}
For the $\mathbb{L}^2$ quasi-tangency, one further requires $\pr{\xi,\eta}\in\Lo{2}{1}\times\Lo{2}{2}$, $\mathbb{P}$-a.s. and the analogous regularity for the correction terms $a,b$. 

\end{proposition}
\begin{proof}[Sketch of the proof]
It is rather clear that the elements in the third point stand for $\pr{\theta_1,\theta_2}$ appearing in Definition \ref{DefLambdaQT}. So, starting from $p^1:=\theta_1-\XX(t+\varepsilon)$, one sets $a^1:=\frac{1}{\varepsilon^{1-\lambda}}\ef{p^1}{t},\ b^1:=\frac{1}{\sqrt{\varepsilon}}\pr{p^1-\ef{p^1}{t}}$. The argument for the $\mathbb{Y}$-component is quite similar and the reasoning gives an equivalence between the sequential characterization and Definition \ref{DefLambdaQT}.
\end{proof}
\begin{remark}
\begin{enumerate}
\item The reader is referred to \cite[Definition 1.1]{AubinDaPrato_90} to see the links of this formulation with the stochastic contingent cones (whenever $\beta_1=\beta_2=0$).
\item Whenever one assumes \eqref{Ass+} to hold true (in addition to \eqref{Ass1} and \eqref{Ass2}), it is easy to see that the fundamental solution $\pr{\XX,\YY}$ shares the $\Lo{2}{1}\times\Lo{2}{2}$ regularity with the initial data $\pr{\xi,\eta}$. As a consequence, one establishes the $\mathbb{L}^2$ quasi-tangency requirements.
\end{enumerate}
\end{remark}
Form now on, and unless stated otherwise, we only look at the $\mathbb{L}^2$-regular notions (quasi-tangency and near viability). 
\subsection{Approximate solutions. Definition and elementary properties}
We consider the following notion of $\varepsilon$-approximate solution compatible at the discretization times with the state constraints.
\begin{definition}\label{DefEpsSol}
For $0\leq t\leq T$, every initial data $\pr{\xi,\eta}\in\mathbb{K}_t\cap\mathbb{L}^2\pr{\Omega,\mathcal{F}_t,\mathbb{P};\Lo{2}{1}\times\Lo{2}{2}}$ and $\varepsilon>0$, an $\varepsilon$-approximate solution for \eqref{Eq0} is a vector $\pr{\bar T,\tau,u,\pr{\phi_1,\phi_2},\pr{\psi_1,\psi_2},\pr{\mathcal{X},\mathcal{Y}}}$ such that
\begin{enumerate}
\item $t\leq\bar{T}\leq T$;
\item the measurable $\tau:\pp{t,\bar{T}}\rightarrow\pp{t,\bar{T}}$ is non-decreasing, non-anticipating and at most $\varepsilon$-delayed i.e.  
\begin{equation*}
 s-\varepsilon\leq\tau(s)\leq s,\ \forall s\in\pp{t,\bar{T}};
\end{equation*} 
\item the control $u\in\mathcal{U}$;
\item the error-capturing processes 
\begin{equation*}
\pr{\phi_1,\phi_2,\psi_1,\psi_2}:\pp{t,\bar{T}}\rightarrow \Hmo{1}\times \Hmo{2}\times \mathcal{L}_2\pr{\Hmo{1}}\times \mathcal{L}_2\pr{\Hmo{2}}
\end{equation*} 
are predictable, take their values in $\Lo{2}{i}$ $\mathbb{P}$-a.s. , for the $\phi_i$, respectively $\mathcal{L}_2\pr{\Lo{2}{i}}$ for the $\psi_i$ (with $i$ given by the index) and satisfy 
\begin{equation*}
\mathbb{E}\pp{\int_t^{\bar T}\normi{\phi_i(l)}{-1}{i}^2dl}\leq \varepsilon\pr{\bar T-t},\ \mathbb{E}\pp{\int_t^{\bar T}\norm{\psi_i(l)}_{\mathcal{L}_2\pr{\Hmo{i}}}^2dl}\leq \varepsilon\pr{\bar T-t};
\end{equation*}
\item the processes $\mathcal{X}, \mathcal{Y}$ are adapted,  $\Lo{2}{1}\times\Lo{2}{2}$-valued,  and satisfy (in the classical $\Hmo{i}$-sense) \begin{align}\label{EqXYCal}
\begin{cases}
d\mathcal{X}(s)&=\pr{\Delta\beta_1\pr{\mathcal{X}(s)}+f_1\pr{\mathcal{X}(\tau(s)),\mathcal{Y}(\tau(s)),u(s)}+\phi_1(s)}ds\\[4pt] & \ \ \ +\pr{\sigma_1\pr{\mathcal{X}(\tau(s)),\mathcal{Y}(\tau(s)),u(s)}+\psi_1(s)}dW(s),\\[4pt]
d\mathcal{Y}(s)&=\pr{\Delta\beta_2\pr{\mathcal{Y}(s)}+f_2\pr{\mathcal{X}(\tau(s)),\mathcal{Y}(\tau(s)),u(s)}+\phi_2(s)}ds\\[4pt] & \ \ \ +\pr{\sigma_2\pr{\mathcal{X}(\tau(s)),\mathcal{Y}(\tau(s)),u(s)}+\psi_2(s)}dW(s),\ s\geq t,\\[4pt]
\mathcal{X}(t)&=\xi,\ \mathcal{Y}(t)=\eta.
\end{cases}
\end{align}
\item For every $s\in\pp{t,\bar{T}}$, the constraint $\pr{\mathcal{X}(\tau(s)),\mathcal{Y}(\tau(s))}\in K$, $\mathbb{P}$-a.s., $\pr{\mathcal{X}(\bar T),\mathcal{Y}(\bar T)}\in \mathbb{K}_{\bar T}$, $\mathbb{P}$-a.s. and 
\begin{equation*}
\mathbb{E}\pp{\normi{\mathcal{X}(\tau(s))-\mathcal{X}(s)}{-1}{1}^2+\normi{\mathcal{Y}(\tau(s))-\mathcal{Y}(s)}{-1}{2}^2}\leq \varepsilon,\ \forall s\in\pp{t,\bar T}.
\end{equation*}
\end{enumerate}
\end{definition}
\begin{proposition}\label{PropPropCal}
We assume \eqref{Ass1} and \eqref{Ass2} to hold true.  Then, there exists a constant $C>0$ such that, for $t\leq s\leq \bar T$ and $\pr{\xi,\eta}\in\mathbb{L}^2\pr{\Omega,\mathcal{F}_t,\mathbb{P};K\cap\pr{\Lo{2}{1}\times\Lo{2}{2}}}$, one has the following properties of an $\varepsilon$-approximate solution as described in Definition \ref{DefEpsSol}.
\begin{enumerate}
\item The following estimate holds true \begin{equation}\label{Estim1cal}\begin{split}
&\mathbb{E}\pp{\sup_{t\leq r\leq s}\normi{\mathcal{X}(r)}{-1}{1}^2+\sup_{t\leq r\leq s}\normi{\mathcal{Y}(r)}{-1}{2}^2+\int_t^s\pr{\norm{\beta_1\pr{\mathcal{X}(l)}}_{\Lo{2}{1}}^2+\norm{\beta_2\pr{\mathcal{Y}(l)}}_{\Lo{2}{2}}^2}dl}\\[4pt]
&\leq C \left( 1 +\mathbb{E}\pp{\normi{\xi}{-1}{1}^2+\int_t^s\normi{\phi_1(l)}{-1}{1}^2dl+\int_t^s\norm{\psi_1(l)}_{\mathcal{L}_2\pr{\Hmo{1}}}^2dl}\right.\\[4pt]
&\left.+\mathbb{E}\pp{\normi{\eta}{-1}{2}^2+\int_t^s\normi{\phi_2(l)}{-1}{2}^2dl+\int_t^s\norm{\psi_2(l)}_{\mathcal{L}_2\pr{\Hmo{2}}}^2dl} \right).
\end{split}\end{equation}
\item Furthermore, if \eqref{Ass+} holds true, then 
\begin{enumerate}
\item $\beta_1$ and $\beta_2$ can be suppressed in the previous inequality such that $\pr{\mathcal{X},\mathcal{Y}}$ has a modification belonging to $\Lo{2}{1}\times\Lo{2}{2}$.
\item This modification (still denoted by $\pr{\mathcal{X},\mathcal{Y}}$) is continuous as a time function with values in $\mathbb{L}^2\pr{\Omega,\mathcal{F},\mathbb{P};\Hmo{1}\times\Hmo{2}}$.
\item The following dependency of the initial data is valid, for $t\leq s\leq \bar T$ \begin{equation}
\label{Estim2cal}
\begin{split}
&\mathbb{E}\pp{\normi{\mathcal{X}(s)-\xi}{-1}{1}^2+\normi{\mathcal{Y}(s)-\eta}{-1}{2}^2}\\[4pt]
&\leq C\pr{1+\mathbb{E}\pp{\norm{\xi}_{\Lo{2}{1}}^2+\norm{\eta}_{\Lo{2}{1}}^2+\int_t^s\pr{\normi{\phi_1(l)}{-1}{1}^2+\normi{\phi_2(l)}{-1}{2}^2}dl}}(s-t)\\[4pt]&+C\mathbb{E}\pp{\int_t^s\pr{\norm{\psi_1(l)}_{\mathcal{L}_2\pr{\Hmo{1}}}^2+\norm{\psi_2(l)}_{\mathcal{L}_2\pr{\Hmo{2}}}^2}dl}
\end{split}
\end{equation}
\end{enumerate}
\end{enumerate}
\end{proposition}
{The proof is very close to the one in Proposition \ref{Prop1}. For our readers' sake, we hint the main supplementary items in the Appendix.}
\subsection{Existence of global approximate solutions}
\begin{theorem}\label{ThAppSol}
Let us assume that \eqref{Ass1}, \eqref{Ass2} and \eqref{Ass+} hold true. Furthermore, we assume that $K\subset\Hmo{1}\times\Hmo{2}$ is a closed set enjoying the quasi-tangency condition with $\lambda=0$ and in the $\mathbb{L}^2$ setting on some interval $\pp{0,T}$. Then, for every initial time $t\in\left[0,T\right)$, for every initial data $\pr{\xi,\eta}\in\mathbb{K}_t$ (and $\Lo{2}{1}\times\Lo{2}{2}$-valued), every time horizon $\tilde{T}\in\pp{t,T}$ and every $\varepsilon\in\pr{0,1}$, there exists a global\footnote{i.e. with the time horizon component $\bar T=\tilde T$} $\varepsilon$-approximate solution in the sense of Definition \ref{DefEpsSol}.
\end{theorem}
\begin{proof}
\textbf{1. From quasi-tangency to local $\varepsilon$-approximate solutions.} \\
We fix $\varepsilon>0$ and pick $0<\varepsilon'<\varepsilon$ (to be made precise at the end of this step).
The quasi-tangency condition at the initial data $(t,\xi,\eta)$ yields the existence of some admissible $u$ (hence the condition 3.  in Definition \ref{DefEpsSol}, $\delta<\varepsilon'$ and
\begin{equation*}
\pr{p^1,p^2}\in \mathbb{L}^2\pr{\Omega,\mathcal{F}_{t+\delta},\mathbb{P};\pr{\Lo{2}{1}\times\Lo{2}{2}}}
\end{equation*}
such that \begin{align}\label{est1_0}\begin{cases}\mathbb{E}\pp{\normi{p^1}{-1}{1}^2+\normi{p^2}{-1}{2}^2}+\frac{1}{\delta}\mathbb{E}\pp{\normi{\ef{p^1}{t}}{-1}{1}^2+\normi{\ef{p^2}{t}}{-1}{2}^2}\leq \varepsilon',\\[8pt]
\pr{\XX\pr{t+\delta}+\sqrt{\delta}p^1,\YY\pr{t+\delta}+\sqrt{\delta}p^2}\in \mathbb{K}_{t+\delta}.
\end{cases}\end{align}
The martingale representation theorem for $\frac{1}{\sqrt{\delta}}p^i$ leads to \[\frac{1}{\sqrt{\delta}}p^i=\phi_i(t)+\int_t^{t+\delta}\frac{1}{\delta}\psi_i(s)dW(s),\]where $\phi_i(t)$ is obtained as the conditional expectation and extended to a function on $\pp{t,t+\delta}$ by setting $\phi_i(s)=\phi_i(t)$.  It follows that 
\begin{equation}\label{est1_1}
\begin{split}
&\mathbb{E}\pp{\int_t^{t+\delta}\pr{\normi{\phi_1(l)}{-1}{1}^2+\normi{\phi_2(l)}{-1}{2}^2+\norm{\psi_1(l)}_{\mathcal{L}_2\pr{\Hmo{1}}}^2+\norm{\psi_2(l)}_{\mathcal{L}_2\pr{\Hmo{2}}}^2}dl}\\[6pt]
&\leq \delta\varepsilon'.
\end{split}
\end{equation}In particular, the condition 4.  in Definition \ref{DefEpsSol} holds true. 
 Furthermore, one sets $\tau(s)=t, \ \forall s\in\pp{t,t+\delta}$ such that $\pr{\mathcal{X}(s),\mathcal{Y}(s)}=\pr{\xi,\eta}\in K$, $\mathbb{P}$-a.s.  The second assertion in \eqref{est1_0} and the construction of $\phi$ and $\psi$ guarantee that \[\pr{\mathcal{X}(t+\delta),\mathcal{Y}(t+\delta)}\in\mathbb{K}_{t+\delta}.\] Finally, owing to \eqref{Estim2cal},  then \eqref{est1_1}, one has
 \begin{align*}
 &\mathbb{E}\pp{\normi{\mathcal{X}(s)-\mathcal{X}(\tau(s))}{-1}{1}^2+\normi{\mathcal{Y}(s)-\mathcal{Y}(\tau(s))}{-1}{2}^2}\\[4pt]&\leq C\pr{1+\mathbb{E}\pp{\norm{\xi}_{\Lo{2}{1}}^2+\norm{\eta}_{\Lo{2}{1}}^2+\int_t^{t+\delta}\pr{\normi{\phi_1(l)}{-1}{1}^2+\normi{\phi_2(l)}{-1}{2}^2}dl}}\delta\\&+C\mathbb{E}\pp{\int_t^{t+\delta}\pr{\norm{\psi_1(l)}_{\mathcal{L}_2\pr{\Hmo{1}}}^2+\norm{\psi_2(l)}_{\mathcal{L}_2\pr{\Hmo{2}}}^2}dl}\\[4pt]
 &\leq C\pr{1+\mathbb{E}\pp{\norm{\xi}_{\Lo{2}{1}}^2+\norm{\eta}_{\Lo{2}{1}}^2}+\varepsilon'}\delta\leq C\pr{1+\mathbb{E}\pp{\norm{\xi}_{\Lo{2}{1}}^2+\norm{\eta}_{\Lo{2}{1}}^2}}\varepsilon'\leq \varepsilon,
 \end{align*}
 explaining the judicious choice of $\varepsilon'$.\\
 
 \textbf{2. Asymptotic behaviour of monotone sequences of approximate solutions.} \\
 As it has been proven in the previous step, the family of $\varepsilon$-approximate solutions is non-empty. This family will be denoted by $\mathcal{A}$ and endowed with a partial order relation 
 \[\pr{\bar T^1,\tau^1,u^1,\pr{\phi_1^1,\phi_2^1},\pr{\psi_1^1,\psi_2^1},\pr{\mathcal{X}^1,\mathcal{Y}^1}}\precsim\pr{\bar T^2,\tau^2,u^2,\pr{\phi_1^2,\phi_2^2},\pr{\psi_1^2,\psi_2^2},\pr{\mathcal{X}^2,\mathcal{Y}^2}},\]if $t\leq\bar T^1\leq\bar T^2$, $u^1=u^2, \ \tau^1=\tau^2,\ \phi_i^1=\phi_i^2,\ \psi_i^1=\psi_i^2$ on $\pp{t, \bar T^1}\times\Omega$ and up to an evanescent set to render indistinguishable the restrictions. We aim at proving that every increasing sequence (indexed by a superscript $n\geq 1$) in $\mathcal{A}$ admits a maximum belonging to $\mathcal{A}$. \\
 
To this purpose, one naturally defines $\bar T:=\underset{n\geq 1}{\sup}\bar T^n$. The stationary case is obvious and will be excluded from the argument. We focus on the framework $\bar T^n<\bar T$. One naturally extends \begin{align*}
\tau(s):=\begin{cases}\tau^n(s),&\textnormal{ if }s\in\pp{t,\bar T^n};\\
\sup_{n\geq 1}\tau^n\pr{\bar T^n},&\textnormal{ if }s=\bar T,\end{cases} 
 \end{align*}keeping the properties in 2. in Definition \ref{DefEpsSol}.\\
 Similarly, one extends $u^n$ by picking some $u_0\in U$ and setting $u(\bar{T})=u_0$. The correction terms $\pr{\phi_1^1,\phi_2^1},\pr{\psi_1^1,\psi_2^1}$ are extended by setting them to $0$ at $s=\bar T$.  Predictability is obvious by definition and the inequality in item 4. in Definition \ref{DefEpsSol} are guaranteed by Fatou's lemma. The fact that $\pr{\mathcal{X},\mathcal{Y}}$ extends $\pr{\mathcal{X}^n,\mathcal{Y}^n}$ is a mere consequence of the uniqueness in \eqref{EqXYCal}.  The continuity of $\mathcal{X},\mathcal{Y}$ (see last item in Proposition \ref{PropPropCal}) and the convergence $\underset{n\rightarrow\infty}{\lim}\tau\pr{\bar T^n}=\tau(\bar T)$ and $\underset{n\rightarrow\infty}{\lim}\bar T^n=\bar T$, together with the closedness of $K$ show that $K\ni\pr{\tilde{X}\pr{\bar T^n},\tilde{Y}\pr{\bar T^n}}$ and $K\ni\pr{\tilde{X}\pr{\tau\pr{\bar T^n}},\tilde{Y}\pr{\tau\pr{\bar T^n}}}$. The same continuity in $\mathbb{L}^2\pr{\Omega, \mathcal{F},\mathbb{P};\Hmo{1}\times\Hmo{2}}$ allows one to pass to the limit as $n\rightarrow\infty$ in the upper-estimate 
\begin{equation*}
\mathbb{E}\pp{\normi{\mathcal{X}(\tau^n(s))-\mathcal{X}(s)}{-1}{1}^2+\normi{\mathcal{Y}(\tau^n(s))-\mathcal{Y}(s)}{-1}{2}^2}\leq \varepsilon,
\end{equation*} to complete the proof. \\

\textbf{3. Existence of a global approximate solution.} \\
We introduce the function $\n{\cdot}:\mathcal{A}\rightarrow\mathbb{R}_+$ defined by $\n{\pr{\bar T,\tau,u,\pr{\phi_1,\phi_2},\pr{\psi_1,\psi_2},\pr{\mathcal{X},\mathcal{Y}}}}:=\bar T.$ One applies Brézis-Browder Theorem (see, for instance \cite[Theorem 2.1.1]{CarjaNeculaVrabie2007}) to deduce the existence of a $\n{\cdot}$-maximal element of $\mathcal{A}$, denoted by $\pr{\bar T^*,\tau^*,u^*,\pr{\phi_1^*,\phi_2^*},\pr{\psi_1^*,\psi_2^*},\pr{\mathcal{X}^*,\mathcal{Y}^*}}$.  Either $\bar T^*=\tilde{T}$ and our claim is proven, or $\bar T^*<\tilde{T}$. In this later case, we apply the first step starting at $\pr{\xi',\eta'}:=\pr{\mathcal{X}\pr{\bar T^*},\mathcal{Y}\pr{\bar T^*}}$ to construct a $\varepsilon$-approximate solution on some non-trivial interval $\pp{\bar T^*,\bar T^*+\delta}$. We emphasize that the linearity in time of condition 4.  in Definition \ref{DefEpsSol} intervenes at this point for the concatenated error processes. It follows that $\bar T^*=\tilde{T}$ (otherwise, the maximality has been contradicted).
\end{proof}
\subsection{The sufficiency result}
\begin{theorem}\label{ThSuf}
Let us assume \eqref{Ass1}, \eqref{Ass2} and \eqref{Ass+} to hold true. Then, if the closed set $K$ is quasi-tangent with $\lambda=0$ for every $t\in\left[0,T\right)$ and every initial data
$\pr{\xi,\eta}$ belonging to the space $\mathbb{K}_t\cap\mathbb{L}^2\pr{\Omega,\mathcal{F}_t,\mathbb{P};\Lo{2}{1}\times\Lo{2}{2}}$, then $K$ is ($\mathbb{L}^2$)-near viable with respect to \eqref{Eq0}.
\end{theorem}
\begin{proof}
We fix $t\in\left[0,T\right)$, $\pr{\xi,\eta}\in\mathbb{K}_t$ with $\Lo{2}{1}\times \Lo{2}{2}$ values and, for the time being, $\varepsilon\in (0,1)$. Theorem \ref{ThAppSol} provides a global $\varepsilon$-approximate solution on $\pp{t,T}$, solution denoted by 
\[\pr{T,\tau,u,\pr{\phi_1,\phi_2},\pr{\psi_1,\psi_2},\pr{\mathcal{X},\mathcal{Y}}}. \] One easily notes (owing to item 6. in Definition \ref{DefEpsSol}) that \begin{equation}\label{est_suff_1}\begin{split}
&d^2\pr{\pr{\X(s),\Y(s)},\mathbb{K}_s}\\[4pt]
&\leq \mathbb{E}\pp{\normi{\X(s)-\mathcal{X}\pr{\tau(s))}}{-1}{1}^2+\normi{\Y(s)-\mathcal{Y}\pr{\tau(s))}}{-1}{2}^2}\\[4pt]
&\leq 2\pr{\mathbb{E}\pp{\normi{\X(s)-\mathcal{X}\pr{s)}}{-1}{1}^2+\normi{\Y(s)-\mathcal{Y}\pr{s)}}{-1}{2}^2}}\\[4pt]& \ \ \ +2\pr{\mathbb{E}\pp{\normi{\mathcal{X}(s)-\mathcal{X}\pr{\tau(s))}}{-1}{1}^2+\normi{\mathcal{Y}(s)-\mathcal{Y}\pr{\tau(s))}}{-1}{2}^2}}\\[4pt]
&\leq 2\pr{\mathbb{E}\pp{\normi{\X(s)-\mathcal{X}\pr{s)}}{-1}{1}^2+\normi{\Y(s)-\mathcal{Y}\pr{s)}}{-1}{2}^2}} +2\varepsilon.
\end{split}\end{equation}For the remaining terms, one begins with writing Itô's formula for $\normi{\X-\mathcal{X}}{-1}{1}^2$ and reasons as in Proposition \ref{Prop2}. Let us briefly explain the treatment of these terms.
We consider $t\leq r\leq s\leq T$. The reader will note that, due to Assumption \eqref{Ass1},
\begin{align*}
&\normi{\X(r)-\mathcal{X}(r)}{-1}{1}^2\\[4pt]
=&2\int_t^r\scal{\Delta\beta_1\pr{X(l)}-\Delta\beta_1\pr{\mathcal X(l)},X(l)-\mathcal X(l)}_{\Hmo{1}}dl\\[4pt]
&+2\int_t^r\scal{f_1\pr{X(l),Y(l),u(l)}-f_1\pr{\mathcal{X}(\tau(l)),\mathcal Y(\tau(l)),u(l)},X(l)-\mathcal X(l)}_{\Hmo{1}}dl\\[4pt]
&+2\int_t^r\scal{\phi_1(l),X(l)-\mathcal X(l)}_{\Hmo{1}}dl+\\[4pt]
&+2\int_t^r\scal{X(l)-\mathcal X(l),\pr{\sigma_1\pr{X(l),Y(l),u(l)}-\sigma_1\pr{\mathcal{X}(\tau(l)),\mathcal Y(\tau(l)),u(l)}-\psi_1(l)}dW_l}_{\Hmo{1}}\\[4pt]
&+\int_t^r\norm{\sigma_1\pr{X(l),Y(l),u(l)}-\sigma_1\pr{\mathcal{X}(\tau(l)),\mathcal Y(\tau(l)),u(l)}-\psi_1(l)}_{\mathcal{L}_2\pr{\Hmo{1}}}^2dl,
\end{align*}
where, in order to avoid long expressions, $\pr{X,Y}=\pr{\X,\Y}$.
For each term,  we bound the scalar product with the product of norms and, for the terms presenting $\tau(l)$, we intercalate $l$. Let us illustrate this on the second term \begin{align*}
&\normi{f_1\pr{\X(l),\Y(l),u(l)}-f_1\pr{\mathcal{X}(\tau(l)),\mathcal Y(\tau(l)),u(l)}}{-1}{1}\normi{\X(l)-\mathcal X(l)}{-1}{1}\\[4pt]
\leq &\normi{\X(l)-\mathcal X(l)}{-1}{1}\normi{f_1\pr{\X(l),\Y(l),u(l)}-f_1\pr{\mathcal{X}(l),\mathcal Y(l),u(l)}}{-1}{1}\\[4pt]&
+\normi{\X(l)-\mathcal X(l)}{-1}{1}\normi{f_1\pr{\mathcal{X}(\tau(l)),\mathcal Y(\tau(l)),u(l)}-f_1\pr{\mathcal{X}(l),\mathcal Y(l),u(l)}}{-1}{1}\\[4pt]
\leq &C\pr{\normi{\X(l)-\mathcal X(l)}{-1}{1}^2+\normi{\Y(l)-\mathcal Y(l)}{-1}{2}^2}\\[4pt]&+C\pr{\normi{\mathcal X(\tau(l))-\mathcal X(l)}{-1}{1}^2+\normi{\mathcal Y(\tau(l))-\mathcal Y(l)}{-1}{2}^2}.
\end{align*}
Here, $C$ only depends on the Lipschitz constant of $f_1$.
As a consequence, we get
\begin{equation*}\begin{split}
&\mathbb{E}\pp{\norm{\X(r)-\mathcal X(r)}_{\Hmo{1}}^2}\\[4pt]
&\leq C\mathbb{E}\pp{\int_t^r \pr{\normi{\X(l)-\mathcal X(l)}{-1}{1}^2+\normi{\Y(l)-\mathcal Y(l)}{-1}{2}^2}dl}\\[4pt]
&+C\int_t^r\mathbb{E}\pp{\pr{\normi{\mathcal X(\tau(l))-\mathcal X(l)}{-1}{1}^2+\normi{\mathcal Y(\tau(l))-\mathcal Y(l)}{-1}{2}^2}}dl
\\[4pt]&+C\mathbb{E}\pp{\int_t^r\normi{\phi_1(l)}{-1}{1}^2dl+\int_t^r\norm{\psi_1(l)}_{\mathcal{L}_2\pr{\Hmo{1}}}^2dl}
\\[4pt]&\leq C\mathbb{E}\pp{\int_t^r \pr{\normi{\X(l)-\mathcal X(l)}{-1}{1}^2+\normi{\Y(l)-\mathcal Y(l)}{-1}{2}^2}dl+\varepsilon}.
\end{split}\end{equation*}The last inequality has been obtained using points 4. and 6. in Definition \ref{DefEpsSol}.  Similar arguments can be developed for $\Y-\mathcal Y$. Putting these together, and using Gronwall's inequality, 
\[\mathbb{E}\pp{\norm{\X(s)-\mathcal X(s)}_{\Hmo{1}}^2+\norm{\Y(s)-\mathcal Y(s)}_{\Hmo{2}}^2}\leq C\varepsilon,\ \forall t\leq s\leq T.\]Going back to \eqref{est_suff_1}, it follows that \[\sup_{t\leq s\leq T}d^2\pr{\pr{\X(s),\Y(s)},\mathbb{K}_s}\leq C\varepsilon,\]which allows one to conclude due to the arbitrariness of $\varepsilon>0$.
\end{proof}
}
{
\section{Application to pseudo-stabilization}\label{SecAppl}
A key method in the asymptotic stabilization of (S)PDEs relies on a "splitting" procedure. This consists in distinguishing between the space generated by the first (generally small) eigenvalues of the differential operator (space referred to as \textbf{unstable}) and the orthogonal component. Then, a first step in achieving stabilizability, is to investigate the existence of controls allowing to drive the unstable part close to $0$ in a quantifiable manner. We offer here an application of the theoretical viability result presented before to this first step (pseudo-stabilization).\\

We consider the case in which $\beta_2\equiv 0,\ f_2(x,y,u)=-cy,$ for some $c>0$, $\sigma_2\equiv 0$ and the set \begin{align}\label{KStab}K=\set{(x,y)\in\Hmo{1}\times\mathbb{R}:\ \underset{1\leq k\leq j}{\sum}\scal{x,e^1_k}_{\Hmo{1}}^2\leq y},\end{align}for some integer $j\geq 1$. The viability of such sets implies that the projection on the eigen-space spanned by the first $j\geq 1$ eigen-functions of $\Delta$ (on $\Lo{2}{1}$ and null on the boundary $\partial{\mathcal{O}_1}$) can be controlled to $0$ with a $c$-exponential speed, i.e.  \[\normi{\Pi_j^1\pr{X^{0,x,y,u}(t)}}{-1}{1}^2\leq ye^{-ct},\ \mathbb{P}\times\mathcal{L}eb-\textnormal{a.s. provided that }\normi{\Pi_j^1(x)}{-1}{1}^2\leq y.\]

Before proceeding, let us point out that Proposition \ref{PropSeqQT} can alternatively be written by applying the martingale representations theorem for $b^i_n$ to find predictable process $\bar{b}_n^1\in\mathcal{L}_2\pr{\Hmo{i}},\ \bar{b}_n^2\in \mathcal{L}_2\pr{G;\mathbb{R}},\ \mathbb{P}$-a.s.  such that $\sqrt{\varepsilon_n}b_n^i=\int_t^{t+\varepsilon_n}\bar{b}_n^i(s)dW(s)$. The properties on $b$ are translated as \begin{align}
\label{Estimbbar}
\lim_{n\rightarrow\infty}\frac{1}{\varepsilon_n}\mathbb{E}\pp{\int_t^{t+\varepsilon_n}\norm{\bar{b}_n^i(s)}_{\mathcal{L}_2\pr{\cdot}}^2ds}=0.
\end{align}Here, $\cdot$ stands for the appropriate spaces i.e.  $\Hmo{1}$ for $b_n^1$ and $\mathbb{R}$ for $b_n^2$.
\begin{remark}
The reader may wonder why $\mathbb{R}$ is considered while the theoretical results are set on $\Lo{2}{2}$. The answer is quite obvious as $\Pi_1^2\pr{\Lo{2}{1}}$ provides us with a one-dimensional space that can (and will) be identified with $\mathbb{R}$.
\end{remark}

For readability purposes, we assume at this point that we deal with an uncontrolled equation \eqref{Eq0} and drop the dependency of $u$.  Therefore, the result hereafter is in the spirit of Nagumo's theorem cf. \cite{Nagumo_42} and the modifications needed in the controlled framework will be explained at the end of the section. \\
For the projections, the eigenvalues and functions, we employ the same notations as in Subsection \ref{SubsecNec4}. \\

Furthermore, let us assume that $K$ is quasi-tangent with $\lambda=0$.  Itô's formula applied for the functional $\pr{x,y}\mapsto \normi{\Pi^1_j(x)}{-1}{1}^2-y$ to \[\pr{\theta_1(s),\theta_2(s)}:=\pr{\mathbb{X}^{t,\xi,\eta}(s)+(s-t)a_n^1+\int_t^s\bar{b}_n^1(l)dW(l),\eta e^{-c(s-t)}+(s-t)a_n^2+\int_t^s\bar{b}_n^2(l)dW(l)},\]yields 
\begin{align*}
&\phiK{\theta_1(s)}{\theta_2(s)}\\
=&\phiK{\theta_1(t)}{\theta_2(t)}++2\int_t^s\pr{c\eta e^{-c(r-t)}-a_n^2}dr\\&+2\int_t^s\sum_{1\leq k\leq j}\scal{\theta_1(r),e_k^1}_{\Hmo{1}}\scal{-\lambda_k^1\beta\pr{\mathbb{X}^{t,\xi,\eta}(r)}+f_1\pr{\xi,\eta}+a_n^1,e_k^1}_{\Hmo{1}}dr\\
&+2\int_t^s\pr{\sum_{1\leq k\leq j}\scal{\theta_1(r),e_k^1}_{\Hmo{1}}\scal{e_k^1,\pr{\sigma_1\pr{\xi,\eta}+\bar b_n^1(r)}dW(r)}_{\Hmo{1}}-\bar{b}_n^2(r)dW(r)}\\&+\int_t^s\norm{\Pi_j^1\pr{\sigma_1\pr{\xi,\eta}+\bar b_n^1(r)}}_{\mathcal{L}_2\pr{\Hmo{1}}}^2dr.
\end{align*}
Note that the left-hand term is non-positive if $\pr{\theta_1,\theta_2}\in K,\ \mathbb{P}$-a.s., while the first term on the right side is null on the boundary $\partial K$. By taking $s=t+\varepsilon_n$ like in Proposition \ref{PropSeqQT}, with initial data belonging to the boundary $\pr{\xi,\eta}\in\partial K,\ \mathbb{P}$-a.s., and relying on the properties of $a_n, b_n$, it follows that
\begin{equation*}\begin{split}
0\geq  &2\int_t^{t+\varepsilon_n}\sum_{1\leq k\leq j}\scal{\mathbb{X}^{t,\xi,\eta}(r),e_k^1}_{\Hmo{1}}\scal{-\lambda_k^1\beta\pr{\mathbb{X}^{t,\xi,\eta}(r)}+f_1\pr{\xi,\eta},e_k^1}_{\Hmo{1}}dr\\&+2c\eta\varepsilon_n+\int_t^{t+\varepsilon_n}\norm{\Pi_j^1\pr{\sigma_1\pr{\xi,\eta}}}_{\mathcal{L}_2\pr{\Hmo{1}}}^2dr+o\pr{\varepsilon_n}\\
&+2\int_t^{t+\varepsilon_n}\sum_{1\leq k\leq j}\scal{\mathbb{X}^{t,\xi,\eta}(r),e_k^1}_{\Hmo{1}}\scal{e_k^1,\sigma_1\pr{\xi,\eta}dW(r)}_{\Hmo{1}}\\
&+2\int_t^{t+\varepsilon_n}\pr{\sum_{1\leq k\leq j}\scal{\mathbb{X}^{t,\xi,\eta}(r),e_k^1}_{\Hmo{1}}\scal{e_k^1,\bar b_n^1(r)}dW(r)_{\Hmo{1}}-\bar{b}_n^2(r)dW(r)},
\end{split}\end{equation*}
where $o$ stands for the classical Landau notation. We now consider $\xi=x\in\Lo{2}{1},\ \eta=y$ on the boundary of $K$ i.e. $y=\normi{\Pi_j^1(x)}{-1}{1}^2$.
\begin{enumerate}
\item By dividing by $\sqrt{\varepsilon_n}$ and using Proposition \ref{Prop1} (assertion 2) and the properties of the Brownian motion,  the associated coefficient must be null i.e.  \begin{equation}\label{Cond2}\sigma_1^*\pr{x,\normi{\Pi_j^1(x)}{-1}{1}^2}\Pi_j^1(x)=0.\end{equation}Please note that the $\bar b$ terms are disappearing due to \eqref{Estimbbar}.
\item By taking the expectancy,  and using, again, Proposition \ref{Prop1}, one gets \begin{equation}\label{interm}\begin{split}
0\geq  &2\varepsilon_n\sum_{1\leq k\leq j}\scal{x,e_k^1}_{\Hmo{1}}\scal{-\lambda_k^1\beta\pr{x}+f_1\pr{x,y},e_k^1}_{\Hmo{1}}\\[4pt]&+2\varepsilon_ncy+\varepsilon_n\norm{\Pi_j^1\pr{\sigma_1\pr{x,y}}}_{\mathcal{L}_2\pr{\Hmo{1}}}^2+o(\varepsilon_n)\\[4pt]
&+2\sum_{1\leq k\leq j}\scal{x,e_k^1}_{\Hmo{1}}\int_t^{t+\varepsilon_n}\mathbb{E}\pp{\scal{-\lambda_k^1\pp{\beta\pr{\mathbb{X}^{t,x,y}(r)}-\beta\pr{x}},e_k^1}_{\Hmo{1}}}dr.
\end{split}\end{equation}
One easily notes that, for some generic constant $C>0$,
\begin{equation*} 
\abs{\scal{\lambda_k^1\pr{\beta\pr{\mathbb{X}^{t,x,y}(r)}-\beta(x)},e_k^1}_{\Hmo{1}}}\leq C\norm{\mathbb{X}^{t,x,y}(r)-x}_{\Lo{2}{1}},
\end{equation*}
and, thus,  
\begin{equation*}
r\mapsto \mathbb{E}\pp{\scal{-\lambda_k^1\beta\pr{\mathbb{X}^{t,x,y}(r))-\beta\pr{x}},e_k^1}_{\Hmo{1}}}
\end{equation*}
is integrable on $\pp{t,t+\varepsilon_0}$. Should $t$ be a Lebesgue point, by dividing \eqref{interm} by $\varepsilon_n$ and allowing $n\rightarrow \infty$,  and recalling that $y=\normi{\Pi_j^1(x)}{-1}{1}^2$, we get the second condition
\begin{equation}
\label{Cond1}\begin{split}
0\geq  &2\sum_{1\leq k\leq j}\scal{x,e_k^1}_{\Hmo{1}}\scal{-\lambda_k^1\beta\pr{x}+f_1\pr{x,\normi{\Pi_j^1(x)}{-1}{1}^2},e_k^1}_{\Hmo{1}}\\[4pt]&+2c\normi{\Pi_j^1(x)}{-1}{1}^2+\norm{\Pi_j^1\pr{\sigma_1\pr{x,\normi{\Pi_j^1(x)}{-1}{1}^2}}}_{\mathcal{L}_2\pr{\Hmo{1}}}^2.
\end{split}
\end{equation}
\end{enumerate}
Let us now prove that \eqref{Cond1} and \eqref{Cond2} are actual sufficient conditions.  To this purpose, let us consider \begin{enumerate}
\item the operator $\Pi_j^\perp:=I-\Pi_j^1$, the projector onto the orthogonal space (in $\Hmo{1}$);
\item For a predictable (fixed) process $\zeta$ taking its values in $\Pi_j^\perp\pr{\Hmo{1}}$, we set \begin{equation}
\begin{cases}
f^{\zeta}\pr{\omega, t, x,y}&:=f_1\pr{x+\zeta(t,\omega),y},\ \sigma^\zeta\pr{\omega, t, x,y}:=\sigma_1\pr{x+\zeta(t,\omega),y},\\[4pt]
F^{\zeta}\pr{\omega, t, x,y}&:=\Pi_j^1\pr{\Delta\beta\pr{x+\zeta(t,\omega)}+f_1\pr{x+\zeta(t,\omega),y}}\\[4pt]
&=\underset{1\leq k\leq j}{\sum}\scal{-\lambda_k^1\beta\pr{x+\zeta(t,\omega)}+f_1\pr{x+\zeta(x,\omega),y},e_k^1}_{\Hmo{1}}e_k^1;\\[4pt]
\Sigma^{\zeta}\pr{\omega, t, x,y}&:=\Pi_j^1\pr{\sigma^\zeta\pr{\omega, t, x,y}};
\end{cases}
\end{equation} for all $\pr{x,y}\in\Pi_j^1\pr{\Hmo{1}}\times \mathbb{R}$. Of course, whenever no confusion is at risk, we will drop the dependency on $\omega\in\Omega$.
\item The equation of interest will be set on the finite-dimensional space \[\mathcal{H}\times\mathbb{R}:=\Pi_j^1\pr{\Hmo{1}}\times \mathbb{R}.\]
\item Due to the choice of $\zeta$, conditions \eqref{Cond1} and \eqref{Cond2} give\begin{equation}
\label{Cond}
\begin{cases}
2\scal{F^{\zeta}\pr{t, x,\norm{x}_{\mathcal{H}}^2},x}_{\mathcal{H}}+2c\norm{x}_{\mathcal{H}}^2+\norm{\Sigma^{\zeta}\pr{t, x,\norm{x}_{\mathcal{H}}^2}}_{\mathcal{L}_2\pr{G;\mathcal{H}}}^2\leq 0;\\[4pt]
\pr{\Sigma^{\zeta}\pr{t, x,\norm{x}_{\mathcal{H}}^2}}^*x=0,\ \forall x\in\mathcal{H},\forall t\in\pp{0,T}.
\end{cases}
\end{equation}
\item Then, one considers, on the finite-dimensional space $\mathcal{H}\times\mathbb{R}$,  the standard equation (with random Lipschitz coefficients) \begin{equation}
\label{EqFin}\begin{cases}
&d^\zeta\mathbf{X}^{t,\xi,\eta}(s)=F^{\zeta}\pr{s,^\zeta\mathbf{X}^{t,\xi,\eta}(s),Y^{t,\xi,\eta}(s)}ds+\Sigma^{\zeta}\pr{s,^\zeta\mathbf{X}^{t,\xi,\eta}(s),Y^{t,\xi,\eta}(s)}dW(s), \\[4pt]
&dY^{t,\xi,\eta}(s)=-cY^{t,\xi,\eta}(s),\ s\in\pp{t,T};\\[4pt]
&^\zeta\mathbf{X}^{t,\xi,\eta}(t)=\xi,\ Y^{t,\xi,\eta}(t)=\eta,
\end{cases}\end{equation}where $\pr{\xi,\eta}\in\mathbb{L}^2\pr{\Omega,\mathcal{F}_t,\mathbb{P};\mathcal{H}\times\mathbb{R}}$.
We draw attention on the fact that the only possible problem ($\beta$ is only Lipschitz in $\Lo{2}{1}$) has now vanished since the contribution on $\Pi_j^\perp\Hmo{1}$ are identical.
\item The reader will immediately note that $\Pi_j^1\pr{X^{t,\xi,\eta}}=^{\Pi_j^\perp\pr{X^{t,\xi,\eta}}}\mathbf{X}^{t,\xi,\eta}$ (i.e., we need to complete an equation by setting $\zeta:=\Pi_j^{\perp}\pr{X^{t,\xi,\eta}}$).
\item The set $K$ can be regarded as a closed convex subset of $\mathcal{H}\times\mathbb{R}$ via \[\mathcal{K}:=\set{\pr{x,y}\in\mathcal{H}\times\mathbb{R}:\ \norm{x}_{\mathcal{H}}^2\leq y}\] which is regular and whose outward normal vectors at boundary points $\pr{x,\norm{x}_{\mathcal{H}}^2}$ are proportional to $\pr{2x,-1}$ (and so is the derivative of the distance function $d_{\mathcal{K}}^2$; see \cite[Example 10]{G18} for analogous computations).
Then, for every fixed $\zeta:= \Pi_j^\perp\pr{X^{t,\xi,\eta}}$, one can apply standard arguments for finite-dimensional systems. The reader is referred to \cite[Theorem 1]{BPQR1998} for the characterization of the viability via properties of $d_{\mathcal{K}}^2$ (equivalent to \cite{AubinDaPrato_90}) to deduce the viability of $\mathcal{K}$. This implies the viability of $K$ and our argument is complete.
\end{enumerate}
\begin{remark}
\begin{enumerate}
\item Let us point out that the cited references (\cite{AubinDaPrato_90}, \cite{G18}, \cite{BPQR1998}) only deal with non-random coefficients, but the result can be extended to random coefficients with similar conclusions (see \cite[Theorem 1]{Shi_2021}).
\item {Assume that the conditions \eqref{Cond1} and \eqref{Cond2} hold true \emph{for every $j\geq 1$.} Starting with $\pr{x,y}$ such that $\normi{x}{-1}{1}^2\leq y$, it follows that, for every $j\geq 1$,  one has $\normi{\Pi^1_j(x)}{-1}{1}^2\leq y$. Then, by the invariance of $K$ (given for $j$ fixed), it follows that $\normi{\Pi^1_j\pr{X^{0,x,y}(t)}}{-1}{1}^2\leq ye^{-ct},\ \mathbb{P}$-a.s.  Since this is valid for every $j\geq 1$, it follows that $\normi{X^{0,x,y}(t)}{-1}{1}^2\leq ye^{-ct},\ \mathbb{P}$-a.s.  In other words, we have obtained a stabilization results not only for the finite-dimensional projection, but, actually, for the entire solution.}
\item Whenever the coefficients $f_1, \sigma_1,\ c$ are controlled, one proves in a similar way that the conditions \eqref{Cond1} and \eqref{Cond2} can be formulated as \begin{equation}\label{CondCtrl}
\begin{split}0\geq \inf&\left\{2\sum_{1\leq k\leq j}\scal{x,e_k^1}_{\Hmo{1}}\scal{-\lambda_k^1\beta\pr{x}+f_1\pr{x,\normi{\Pi_j^1(x)}{-1}{1}^2,u},e_k^1}_{\Hmo{1}} \right.
\\[4pt]&+2c(u)\normi{\Pi_j^1(x)}{-1}{1}^2+\norm{\Pi_j^1\pr{\sigma_1\pr{x,\normi{\Pi_j^1(x)}{-1}{1}^2,u}}}_{\mathcal{L}_2\pr{\Hmo{1}}}^2;\\[6pt]
&\left. u\in U\textnormal{ such that }\sigma_1^*\pr{x,\normi{\Pi_j^1(x)}{-1}{1}^2,u}\Pi_j^1(x)=0 \right\}.\end{split}
\end{equation}The reader is referred to \cite[Theorem A1]{BCQ2001} for connections with the finite-dimensional diffusions and a similar qualitative result for that framework.
\item Refinements linked to $c$ with additional dependence on the $x,y$ components can also be envisaged.
\end{enumerate}
\end{remark}
}
\section{Appendix}\label{SecAppendix}
{We provide here some of the technical considerations leading to the inequalities used throughout the paper.}
\subsection{A1. Proofs for Propositions \ref{Prop1}, \ref{PropContInitialData} and \ref{Prop2}}
\begin{proof}[Proof of Proposition \ref{Prop1}]
Let us consider $t\leq r\leq s\leq T$. \\
1. The reader will note that, due to Assumption \eqref{Ass1} (see also \eqref{Estim0}) and using Itô's formula,
\begin{align*}
&\norm{\XX(r)}_{\Hmo{1}}^2\\[4pt] 
=&\norm{\xi}_{\Hmo{1}}^2+2\int_t^r\scal{\Delta\beta_1\pr{\XX(l)},\XX(l)}_{\Hmo{1}}dl\\[4pt]  
&+2\int_t^r\scal{f_1\pr{\xi,\eta,u(l)},\XX(l)}_{\Hmo{1}}dl\\[4pt] 
&+2\int_t^r\scal{\XX(l),\sigma_1\pr{\xi,\eta,u(l)}dW_l}_{\Hmo{1}}+\int_t^r\norm{\sigma_1\pr{\xi,\eta,u(l)}}_{\mathcal{L}_2\pr{\Hmo{1}}}^2dl\\[4pt] 
\leq &2\abs{\int_t^r\scal{\XX(l),\sigma_1\pr{\xi,\eta,u(l)}dW_l}_{\Hmo{1}}}-2\bar \alpha_1\int_t^r\norm{\beta_1\pr{\XX(l)}}_{\Lo{2}{1}}^2dl\\[4pt] 
&+C\pp{\int_t^r\norm{\XX(l)}_{\Hmo{1}}^2dl+\pr{(r-t)\vee 1}\pr{1+\normi{\xi}{-1}{1}^2+\normi{\eta}{-1}{2}^2}}.
\end{align*}Then, owing to Burkholder-Davis-Gundy inequality and using, again,  \eqref{Ass1},  followed by Gronwall's inequality, it follows that\begin{equation*}\begin{split}
&\mathbb{E}\pp{\sup_{t\leq r\leq s}\norm{\XX(r)}_{\Hmo{1}}^2+\int_t^s\norm{\beta_1\pr{\XX(l)}}_{\Lo{2}{1}}^2dl}\\[4pt] 
&\leq C\pr{1+\mathbb{E}\pp{\norm{\xi}_{\Hmo{1}}^2+\norm{\eta}_{\Hmo{2}}^2}}.\end{split}
\end{equation*}(We have implicitly used $\bar\alpha_1>0$. ) A similar argument applies to $\YY$.
Our first claim follows.\\
2.  Before moving to the remaining claims, let us just point out that, whenever $\beta\neq 0$ and \eqref{Ass1} holds true,  one can drop the functions $\beta_1$ and $\beta_2$ in the left-hand member of \eqref{L2estim}. {The case $\beta=0$ is a standard Lipschitz-one and the equation (and the well-posedness) are much simpler. We will only explain the arguments for the more challenging part when $\beta\neq 0$.}\\
To prove the second assertion, we begin with writing down Itô's formula for $\norm{\XX-\xi}_{\Hmo{1}}^2$ and $\norm{\YY-\eta}_{\Hmo{2}}^2$. 
\begin{equation}\label{E1}\begin{split}
&\norm{\XX(r)-\xi}_{\Hmo{1}}^2\\[4pt]  
=&2\int_t^r\scal{\Delta\beta_1\pr{\XX(l)},\XX(l)-\xi}_{\Lo{2}{1}^*,\Lo{2}{1}}dl
\\[4pt]&+2\int_t^r\scal{f_1\pr{\xi,\eta,u(l)},\XX(l)-\xi}_{\Hmo{1}}dl\\[4pt] 
&+2\int_t^r\scal{\XX(l)-\xi,\sigma_1\pr{\xi,\eta,u(l)}dW(l)}_{\Hmo{1}}+\int_t^r\norm{\sigma_1\pr{\xi,\eta,u(l)}}_{\mathcal{L}_2\pr{\Hmo{1}}}^2dl
\\[4pt] \leq &-2\alpha_1\int_t^r\norm{\XX(l)-\xi}_{\Lo{2}{1}}^2dl\\[4pt]
+&C\int_t^r\norm{\XX(l)-\xi}_{\Hmo{1}}^2dl+2\int_t^r\norm{\Delta\beta_1\pr{\xi}}_{\Lo{2}{1}^*}\norm{\XX(l)-\xi}_{\Lo{2}{1}}dl\\[4pt]
&+C\abs{\int_t^r\scal{\XX(l)-\xi,\sigma_1\pr{\xi,\eta,u(l)}dW(l)}_{\Hmo{1}}}+C(r-t)\pr{1+\norm{\xi}_{\Hmo{1}}^2+\norm{\eta}_{\Hmo{2}}^2}\\[4pt] 
\leq &C\int_t^r\norm{\XX(l)-\xi}_{\Hmo{1}}^2dl+C\int_t^r\norm{\Delta\beta_1\pr{\xi}}_{\Lo{2}{1}^*}^2dl\\[4pt] 
&+C\abs{\int_t^r\scal{\XX(l)-\xi,\sigma_1\pr{\xi,\eta,u(l)}dW(l)}_{\Hmo{1}}}+C(r-t)\pr{1+\norm{\xi}_{\Hmo{1}}^2+\norm{\eta}_{\Hmo{2}}^2}.
\end{split}\end{equation}
Although classical, let us explain how to deal with a quadratic variation term (linked to the $W$ term).  Using the assumptions on $\sigma_1$, it follows that 
\begin{equation}\label{E2}\begin{split}
&\mathbb{E}\pp{\pr{\int_t^s\norm{\XX(l)-\xi}_{\Hmo{1}}^2\norm{\sigma_1\pr{\xi,\eta,u(l)}}_{\Hmo{1}}^2dl}^{\frac{1}{2}}}\\[4pt]
&\leq C\mathbb{E}\pp{\sup_{t\leq r\leq s}\norm{\XX(l)-\xi}_{\Hmo{1}}\pr{(s-t)\pr{1+\norm{\xi}_{\Hmo{1}}^2+\norm{\eta}_{\Hmo{2}}^2}}^{\frac{1}{2}}}   \\[4pt] 
&\leq \frac{1}{4}\mathbb{E}\pp{\sup_{t\leq r\leq s}\norm{\XX(l)-\xi}_{\Hmo{1}}^2}+C(s-t)\pr{1+\mathbb{E}\pp{\norm{\xi}_{\Hmo{1}}^2+\norm{\eta}_{\Hmo{2}}^2}}.
\end{split}\end{equation}
Owing to Burkholder-Davis-Gundy inequality and to the arguments in \eqref{E1} and \eqref{E2}, completed with Hölder's inequality, it follows that 
\begin{equation*}\begin{split}
&\mathbb{E}\pp{\sup_{t\leq r\leq s}\norm{\XX(r)-\xi}_{\Hmo{1}}^2} \\[4pt]
\leq \;  & \; C   \mathbb{E}\pp{\int_t^s\norm{\XX(l)-\xi}_{\Hmo{1}}^2dl}  \\[4pt]
& + C(s-t)\pr{1+\mathbb{E}\pp{\norm{\xi}_{\Hmo{1}}^2+\norm{\eta}_{\Hmo{2}}^2+\norm{\Delta\beta_1\pr{\xi}}_{\Lo{2}{1}^*}^2}}  \\[4pt]
\leq \;  & \; C  \mathbb{E}\pp{\int_t^s\norm{\XX(l)-\xi}_{\Hmo{1}}^2dl}  \\[4pt]
& +C(s-t)\pr{1+\mathbb{E}\pp{\norm{\xi}_{\Hmo{1}}^2+\norm{\eta}_{\Hmo{2}}^2+\norm{\xi}_{\Lo{2}{1}}^2}}.
\end{split}
\end{equation*}
The last inequality follows from the last assertion in Remark \ref{Rem1}.  The reader is invited to note that the last term is $0$ if $\beta_1\equiv0$.  Owing to Gronwall's inequality and to the first assertion in our proposition\footnote{Please take a look at the remark opening the proof of this assertion to see why $\beta_1$ is dropped},  if \eqref{Ass+} holds true,  it follows that
\begin{equation*}\begin{split}
& \mathbb{E}\pp{\sup_{t\leq r\leq s}\norm{\XX(r)-\xi}_{\Hmo{1}}^2} \\[6pt]
& \leq C\pr{1+\mathbb{E}\pp{\norm{\xi}_{\Lo{2}{1}}^p+\norm{\xi}_{\Hmo{1}}^2+\norm{\eta}_{\Hmo{2}}^2}}(s-t).
\end{split}
\end{equation*}
The argument on $\norm{\YY-\eta}_{\Hmo{2}}^2$ is quite similar. 
\end{proof}

\begin{proof}[Proof of Proposition \ref{PropContInitialData}]
Itô's formula written for $\normi{\X-\Xp}{-1}{1}^2$ yields \begin{equation*}
\begin{split}
&\normi{\X(r)-\Xp(r)}{-1}{1}^2+2\bar \alpha_1\int_t^r\norm{\beta_1\pr{\X(l)}-\beta_1\pr{\Xp(l)}}_{\Lo{2}{1}}^2dl\\[4pt]
\leq &\normi{\xi-\xi'}{-1}{1}^2\\[4pt]
& \begin{split} +2\int_t^r & \left\langle f_1\pr{\X(l),\Y(l),u(l)}-f_1\pr{\Xp(l),\Yp(l),u(l)},\right. \\[4pt] & \left.   \X(l)-\Xp(l) \right\rangle_{\Hmo{1}} dl \end{split} \\[4pt]
& \begin{split}+2\int_t^r & \left\langle \X(l)-\Xp(l),\right. \\[4pt]& \left.  \sigma_1\pr{\X(l),\Y(l),u(l)}-\sigma_1\pr{\Xp(l),\Yp(l),u(l)}dW(l) \right\rangle_{\Hmo{1}}  \end{split} \\[4pt]
&+\int_t^r\norm{\sigma\pr{\X(l),\Y(l),u(l)}-\sigma_1\pr{\Xp(l),\Yp(l),u(l)}}_{\mathcal{L}_2\pr{\Hmo{1}}}^2dl.
\end{split}\end{equation*}
Using Burkholder-Davis-Gundy inequality, one gets, for every $s\in\pp{t,T}$, 
\begin{align*}&\mathbb{E}
\left[\sup_{t\leq r\leq s}\normi{\X(r)-\Xp(r)}{-1}{1}^2\right.\\[4pt]
&\left.+2\bar \alpha_1\int_t^s\norm{\beta_1\pr{\X(l)}-\beta_1\pr{\Xp(l)}}_{\Lo{2}{1}}^2dl\right] \\[4pt]
&\leq C\mathbb{E}\pp{\int_t^s\pr{\normi{\X(l)-\Xp(l)}{-1}{1}^2+\normi{\Y(r)-\Yp(r)}{-1}{2}^2}dl}\\[4pt]
&+\mathbb{E}\pp{\norm{\xi-\xi'}_{\Hmo{1}}^2}.
\end{align*}
A similar inequality is valid for the term $\normi{\Y-\Yp}{-1}{2}^2$ and the conclusion follows by summing the two and applying Gronwall's inequality. The second inequality is quite similar.
\end{proof}

\begin{proof}[Proof of Proposition \ref{Prop2}]
As in the previous results,  we consider $t\leq r\leq s\leq T$ and we fix an admissible control $u\in\mathcal{U}$.
The reader will note that, due to Assumption \eqref{Ass1},
\begin{align*}
&\normi{\XX(r)-\X(r)}{-1}{1}^2\\[4pt]
=&2\int_t^r\scal{\Delta\beta_1\pr{\XX(l)}-\Delta\beta_1\pr{\X(l)},\XX(l)-\X(l)}_{\Hmo{1}}dl\\[4pt]
&+2\int_t^r\scal{f_1\pr{\xi,\eta,u(l)}-f_1\pr{\X(l),\Y(l),u(l)},\XX(l)-\X(l)}_{\Hmo{1}}dl\\[4pt]
&+2\int_t^r\scal{\XX(l)-\X(l),\sigma_1\pr{\xi,\eta,u(l)}-\sigma_1\pr{\X(l),\Y(l),u(l)}dW_l}_{\Hmo{1}}\\[4pt]
&+\int_t^r\norm{\sigma\pr{\xi,\eta,u(l)}-\sigma_1\pr{\X(l),\Y(l),u(l)}}_{\mathcal{L}_2\pr{\Hmo{1}}}^2dl
\end{align*}
The last term is dealt with via the Lipschitz continuity of $\sigma_1$. For the second term, by dropping the dependency on the time parameter $l$, one writes \begin{align*}
&\normi{f_1\pr{\xi,\eta,u}-f_1\pr{\X,\Y,u}}{-1}{1}\normi{\XX-\X}{-1}{1}\\[4pt]
\leq &\normi{\XX-\X}{-1}{1}\normi{f_1\pr{\xi,\eta,u}-f_1\pr{\XX,\YY,u}}{-1}{1}\\[4pt]&
+\normi{\XX-\X}{-1}{1}\normi{f_1\pr{\XX,\YY,u}-f_1\pr{\X,\Y,u}}{-1}{1}\\[4pt]
\leq &\pp{f_1}_1\pr{\normi{\xi-\XX}{-1}{1}+\normi{\eta-\YY}{-1}{2}}\normi{\XX-\X}{-1}{1}\\[4pt]&+\pp{f_1}_1\pp{\normi{\XX-\X}{-1}{1}^2+\normi{\YY-\Y}{-1}{2}\normi{\XX-\X}{-1}{1}}\\[4pt]
\leq &C\pr{\normi{\XX-\X}{-1}{1}^2+\normi{\YY-\Y}{-1}{2}^2+\normi{\xi-\XX}{-1}{1}^2}\\[4pt] &+C\normi{\eta-\YY}{-1}{2}^2.
\end{align*}
As a consequence, returning to the previous inequality, and recalling that the assumption \eqref{Ass+} holds true, one has 
\begin{align*}
&\normi{\XX(r)-\X(r)}{-1}{1}^2+\int_t^r\norm{\XX(l)-\X(l)}_{\Lo{2}{1}}^{2}dl\\[4pt]
\leq &C\int_t^r\pr{\norm{\XX(l)-\X(l)}^2_{\Hmo{1}}+\normi{\YY(l)-\Y(l)}{-1}{2}^2}dl\\[4pt]&+
\abs{\int_t^r\scal{\XX(l)-\X(l),\pr{\sigma_1\pr{\xi,\eta,u(l)}-\sigma_1\pr{\X(l),\Y(l),u(l)}dW_l}}_{\Hmo{1}}}\\[4pt]
&+C\int_t^r\pr{\norm{\XX(l)-\xi}_{\Hmo{1}}^2+\norm{\YY(l)-\eta}_{\Hmo{2}}^2}dl.
\end{align*}
Proceeding as in the previous proposition (using Burkholder-Davis-Gundy inequality followed by Gronwall's inequality), we get
\begin{equation}\begin{split}
&\mathbb{E}\pp{\sup_{t\leq r\leq s}\norm{\XX(r)-\X(r)}_{\Hmo{1}}^2+\int_t^s\norm{\XX(l)-\X(l)}_{\Lo{2}{1}}^{2}dl}\\[4pt]
\leq \;\;& C\mathbb{E}\int_t^s\left(\norm{\Y(l)-\YY(l)}_{\Hmo{2}}^2+\norm{\XX(l)-\xi}_{\Hmo{1}}^2 \right. \\[4pt]
&\qquad \quad \left.+\norm{\YY(l)-\eta}_{\Hmo{2}}^2 \right] dl.
\end{split}\end{equation}
Similar arguments can be developed for $\YY-\Y$.
The conclusion follows by taking the sum of the estimates on these differences, applying again Gronwall's inequality and due to the last two assertions in Proposition \ref{Prop1}.
\end{proof}

\subsection{A2. Elements of proof for Proposition \ref{PropPropCal}}
\begin{proof}[Sketch of the proof of Proposition \ref{PropPropCal}]
As we have already specified, we only give some hints to the main modifications as this result is quasi-identical in treatment to the one in Proposition \ref{Prop1}.\\
The first assertion follows the same arguments as the ones used in Proposition \ref{Prop1} and so does 2. (a). The particular impact on regularity follows from the integral term containing $\Lo{2}{i}$-norms being bounded.\\
To prove the continuity claim, one proceeds as in Proposition \ref{Prop1}, the last two assertions. Let us fix $t\leq r\leq s\leq \bar T$. Then, is is clear that, on $\pp{r,s}$, the processes $\mathcal{X},\mathcal{Y}$ satisfy \eqref{EqXYCal} with $\Lo{2}{1}\times \Lo{2}{2}$-regular initial data $\mathcal{X}(r),\mathcal{Y}(r)$.  We merely hint to the treatment of the term involving $\Delta$ and refer the reader to the proof of Proposition \ref{Prop1} to complete the elements. For $r\leq l\leq s$,
\begin{align*}
&\int_r^s\scal{\Delta\beta_1\pr{\mathcal{X}(l)},\mathcal{X}(l)-\mathcal{X}(r)}_{\Hmo{1}}dl\\[4pt]
&\leq \int_r^s\norm{\Delta\beta_1\pr{\mathcal{X}(l)}}_{\Lo{2}{1}^*}\norm{\mathcal{X}(l)-\mathcal{X}(r)}_{\Lo{2}{1}}dl\\[4pt]
\leq &C\pr{\int_r^s\norm{\mathcal{X}(l)}_{\Lo{2}{1}}^2dl}^{\frac{1}{2}}\pr{\pr{\int_r^s\norm{\mathcal{X}(l)}_{\Lo{2}{1}}^2dl}^{\frac{1}{2}}+\pr{\norm{\mathcal{X}(r)}_{\Lo{2}{1}}^2}^{\frac{1}{2}}}.
\end{align*}
The conclusion follows by noting that $\underset{\pr{s-r}\rightarrow 0}{\lim}\int_r^s\mathbb{E}\pp{\norm{\mathcal{X}(l)}_{\Lo{2}{1}}^2dl}=0$ due to the integrability on $\pp{t,\bar T}$ of $\mathbb{E}\pp{\norm{\mathcal{X}(l)}_{\Lo{2}{1}}^2dl}$.\\
Finally, for the dependency on the initial data, one writes (as it was already the case in Proposition \ref{Prop1}), Itô's formula for $\normi{\mathcal{X}-\xi}{-1}{1}^2$). The new terms appearing in the right-hand side are dealt with as follows.  First, for $t\leq r\leq s$,
\begin{align*}
&\mathbb{E}\pp{\int_t^r\scal{\phi_1(l),\mathcal{X}(l)-\xi}_{\Hmo{1}}dl}\leq \mathbb{E}\pp{\int_t^s\normi{\phi_1(l)}{-1}{1}dl\times \sup_{t\leq l\leq r}\norm{\mathcal{X}(l)-\xi}_{\Hmo{1}}}\\&\leq \frac{1}{4\delta}\mathbb{E}\pp{\pr{\int_t^s\normi{\phi_1(l)}{-1}{1}dl}^2}+\delta\mathbb{E}\pp{\sup_{t\leq l\leq r}\norm{\mathcal{X}(l)-\xi}_{\Hmo{1}}^2}\\
&\leq (s-t)\frac{1}{4\delta}\mathbb{E}\pp{\int_t^s\normi{\phi_1(l)}{-1}{1}^2dl}+\delta\mathbb{E}\pp{\sup_{t\leq l\leq r}\norm{\mathcal{X}(l)-\xi}_{\Hmo{1}}^2}.
\end{align*}
As usual, $\delta$ is chosen small enough to be compensated by the left-hand side.  Second, the stochastic term involving $\psi_1$ i.e.  $\abs{\int_t^r\scal{\mathcal{X}(l)-\xi,\psi_1(l)dW(l)}_{\Hmo{1}}}$ is dealt with through Burkholder-Davis-Gundy inequality which provides an upper bound linked to \begin{align*}
&\mathbb{E}\pp{\pr{\int_t^s\normi{\mathcal{X}(l)-\xi}{-1}{1}^2\norm{\psi_2(l)}_{\mathcal{L}_2\pr{\Hmo{1}}}^2dl}^{\frac{1}{2}}}\\&\leq \frac{1}{4\delta}\mathbb{E}\pp{\int_t^s\norm{\psi_2(l)}_{\mathcal{L}_2\pr{\Hmo{1}}}^2dl}+\delta\mathbb{E}\pp{\sup_{t\leq l\leq s}\norm{\mathcal{X}(l)-\xi}_{\Hmo{1}}^2},
\end{align*}with $\delta$ chosen as before.
We have now accounted for all the additional contributions in \eqref{Estim2cal} with respect to what appeared in Proposition \ref{Prop1}.
\end{proof}

\section*{Acknowledgements} I.M.  would like to thank the colleagues at the LMI, Normandie University, INSA de Rouen Normandie, for a pleasant stay at their department where part of this work was done.

\bibliographystyle{plain}
\bibliography{Bibl_2021.bib}
\end{document}